\magnification=\magstephalf

\input amstex
\documentstyle{amsppt}
\hoffset 1.25truecm
\hsize=12.4 cm
\vsize=19.7 cm
\TagsOnRight


\font\tenscr=rsfs10 
\font\sevenscr=rsfs7 
\font\fivescr=rsfs5 
\skewchar\tenscr='177 \skewchar\sevenscr='177
\skewchar\fivescr='177
\newfam\scrfam \textfont\scrfam=\tenscr \scriptfont\scrfam=\sevenscr
\scriptscriptfont\scrfam=\fivescr
\def\scr{\fam\scrfam}

\def\az{\alpha}  \def\bz{\beta}
    
    \def\fz{\varphi}
\def\gz{\gamma}  \def\kz{\kappa}
\def\lz{\lambda} 
\def\nz{\nu}

        \def\sz{\sigma}
        \def\uz{\theta}
\def\vz{\varepsilon} 
\def\zz{\zeta}

\def\llz{\Lambda} \def\ddz{\Delta}
    
\def\ggz{\Gamma}

\def\flsh{\flushpar}

\def\qd{\quad}
\def\qqd{\qquad}

\def\lmm{Lemma}
\def\prp{Proposition}
\def\thm{Theorem}

\def\prf{Proof}

\def\d{\text{\rm d}}
\def\var{\text{\rm Var}}
\def\cov{\text{\rm Cov}}
\def\ent{\text{\rm Ent}}
\def\hes{\text{\rm Hess}}
\def\le{\leqslant}
\def\ge{\geqslant}

 \topmatter
\topinsert \captionwidth{10 truecm}
\flsh Acta Math. Sin. New Ser.
\flsh 2008, 24(5), pp. 705--736
\endinsert
\title{Spectral gap and logarithmic Sobolev constant for continuous spin systems}\endtitle
\rightheadtext{Spectral gap and logarithmic Sobolev constant}
\author{Mu-Fa Chen}\endauthor
\affil{(School of Mathematical Sciences, Beijing Normal University, Beijing 100875, P.R. China)\\
        E-mail: mfchen\@bnu.edu.cn\\
        September 19, 2007}\endaffil
\thanks {Research supported in part by the Creative Research Group Fund of the
         National Natural Science Foundation of China (No. 10121101) and by the ``985'' project from the Ministry of Education in China.}\endthanks
\abstract{The aim of
this paper is to study the spectral gap and the logarithmic Sobolev constant for continuous spin
systems. A simple but general result for estimating the spectral gap of
finite dimensional systems is given by Theorem 1.1, in terms of the spectral gap
for one-dimensional marginals. The study of the topic provides us a
chance, and it is indeed another aim of the paper, to justify the
power of the results obtained previously. The exact order in dimension one (Proposition 1.4), and then the precise leading order
and the explicit positive regions of the spectral gap and the logarithmic Sobolev constant
 for two typical infinite-dimensional models are
presented (Theorems 6.2 and 6.3). Since we are interested in explicit estimates, the computations
become quite involved. A long section (Section 4) is devoted to the study of the spectral gap
in dimension one.}\endabstract
\endtopmatter

\document

\subhead{1. Introduction}
\endsubhead
The local Poincar\'e inequalities (equivalently, spectral gaps) and
logarithmic Sobolev inequalities for unbounded continuous spin
systems have recently obtained a lot of attention by many authors
[1]--[11]. For the present status of the study and further
references, the readers may refer to the comprehensive survey
article [7]. In the most of the publications, the authors consider
mainly the perturbation regime with convex phase at infinity. More
recently, the non-convex phase is treated for a class of spin
systems based on a criterion for the weighted Hardy inequalities.

The main purpose of this paper is to propose a general formula for
the local spectral gaps of continuous spin systems. Let us start
from finite dimensions. Let $U\in C^\infty  (\Bbb R^n) $ satisfy
$Z:=\int_{\Bbb R^n} e^{- U} \d x< \infty$ and set $\d
\mu_U^{}=e^{-U}\d x/Z $. Throughout the paper, we use a particular
notation $x_{\backslash i}:=(x_1, \ldots, x_{i-1},x_{i+1},\ldots,
x_n)\in \Bbb R^{n-1} $, obtained from $x:=(x_1, x_2, \ldots,
x_n)\in \Bbb R^n$ by removing the $i$th component. Clearly, the
conditional distribution of $ x_i$ given $x_{\backslash i}$ under $\mu_U^{}$ is as follows:
$$ \mu_U^{x_{\backslash i}} (\d x_i)=e^{-U}\d x_i/ Z(x_{\backslash i}), \tag 1.1$$
where $ Z(x_{\backslash i})=\int_{\Bbb R} e^{-U(x)} \d x_i$. The
measure $\mu_U^{x_{\backslash i}}$ is the invariant probability
measure of the one-dimensional diffusion process, corresponding to
the operator $ L_i^{x_{\backslash i}}=\d^2/\d x_i^2 - \partial_i
U\, \d/\d x_i$.

Let $L=\Delta -\langle \nabla U, \nabla\rangle$. Recall that the
spectral gap $ \lambda_1 (L)=\lambda_1 (U)$ is the largest
constant $\kappa $ in the following Poincar\'e inequality
$$ \kappa \var_{\mu_U^{}} (f)\le
\int_{\Bbb R^n} |\nabla f|^2 \d \mu_U^{} =: D(f), \qquad f\in
C_0^\infty (\Bbb R^n),\tag 1.2$$
where $\var_{\mu_U^{}} (f)$ is the
variation of $f$ with respect to $\mu_U^{} $ and $ C_0^\infty
(\Bbb R^n)$ is the set of smooth functions with compact supports.

Denote by $\lambda_1^{x_{\backslash i}}= \lambda_1\big(
L_i^{x_{\backslash i}}\big)$ the spectral gap of the
one-dimensional operator $ L_i^{x_{\backslash i}}$:
$$ \lambda_1^{x_{\backslash i}} \var_{\mu_U^{ x_{\backslash i}}} (f)\le
\int_{\Bbb R} {f'}^2 \d \mu_U^{ x_{\backslash i}}, \qquad f\in
C_0^\infty (\Bbb R).   \tag 1.3$$ Then, we can state our
variational formula for the lower bounds of $\lambda_1(U)$ as
follows.

\proclaim{Theorem 1.1} Define
$$
\big({\widetilde \hes}(U)\big)_{ij}=
\cases \lambda_1^{x_{\backslash i}}, \quad &\text{if }i=j \\
\partial_{ij}U, \quad &\text{if }i\ne j,
\endcases$$
where $({\hes}(U))_{ij}=\partial_{ij}U:=\partial^2
U/\partial_{x_i}\partial_{x_j}$. Then we have
$$\align
\lambda_1 (U)&\ge \inf_{x\in \Bbb R^n} \lambda_{\min} \big({\widetilde \hes}(U)(x)\big)\\
&\ge \inf_{x\in \Bbb R^n} \sup_{w} \min_{1\le i\le n}
\bigg(\lambda_1^{x_{\backslash i}}- \sum_{j: j\ne i} |\partial_{ij}
U (x)| w_j/w_i \bigg), \tag 1.4 \endalign$$ where $w=(w_i)_{i=1}^n$
varies over all positive sequences.
\endproclaim

Setting $ w_i \equiv 1$ in (1.4), it follows that
$$\lambda_1(U)
\ge \inf_{x\in \Bbb R^n} \min_{1\le i\le n}
\bigg(\lambda_1^{x_{\backslash i}} - \sum_{j: j\ne i} |\partial_{ij}
U(x)|\bigg) \ge  \min_{1\le i\le n} \bigg[\inf_{x\in \Bbb
R^n}\lambda_1^{x_{\backslash i}} - \sum_{j: j\ne i} \|\partial_{ij}
U\|_{\infty }\bigg].
$$
The last lower bound is more or less the estimate given in [5] and
[7], goes back to [3].

The supremum over $w$ in (1.4) comes from a variational formula
for the principal eigenvalue of a symmetric $Q$-matrix (cf. \S3
for more details). The use of the variational formula is
necessary, since the principal eigenvalue is not computable in
general for a large scale matrix.

The essential point for which (1.4) is valuable is that we now have
quite complete knowledge about the spectral gap in dimensional one.
For instance, as a consequence of part (1) of Theorem 3.1 in [12],
we have
$$\lambda_1^{x_{\backslash i}}\ge
\sup_f \inf_{x_i\in \Bbb R}
\bigg\{\partial_{ii}U (x)-\frac{f''(x_i)-
\partial_i U(x)f'(x_i)}{f(x_i)}\bigg\},\tag 1.5$$
where $f$ varies over all positive functions in $C^2 (\Bbb
R)$. In particular, setting $f=1$, we get
$$\lambda_1^{x_{\backslash i}}\ge
\inf_{x_i\in \Bbb R} \partial_{ii}U (x). \tag 1.6 $$ When
$\partial_{ii}U (x)=u''(x_i)$ for some $u\in C^2 (\Bbb R)$,
independent of $i$, (1.6) leads to the so-called convex phase
condition ``$\inf_{x\in \Bbb R} u''(x)>0 $.'' Since a local
modification of $u$ does not change the positiveness of $\lambda_1$,
the convex condition can be replaced by $ \varliminf_{|x|\to\infty }
u''(x)>0$ (i.e., the convexity at infinity) as proved in [12;
Corollary 3.5], see also Theorem 4.1 below. However, the last
condition is still not necessary as shown by [12; Example 3.11 (3):
$u'(x)= \gamma x (\gamma +\cos x)^{-1}$ for some $\gamma>1$] and [5;
Proposition 4.4] (see also Example 2.5 below). A more careful
examination of spectral gap in dimension one is delayed to \S 4.

It is possible to avoid the use of test functions $w$ and $f$ in
(1.4) and (1.5), respectively. To see this, we introduce
an explicit lower estimate of $\lambda_1 (U) $. For this, we need
additional notations. Choose a practical $\eta_i^{x_{\backslash i}}\le \lambda_1^{x_{\backslash i}}$, as bigger as possible, and define
$$\align
&{\matrix
s_i(x)=\eta_i^{x_{\backslash i}}- \sum_{j: j\ne i}|\partial_{ij}
U(x)|,
\qquad \qqd &\;\;\;{\underline s} (x)=\min_{1\le i\le n} s_i(x),\\
q_i(x)=\eta_i^{x_{\backslash i}}-{\underline s}(x) , \qquad  &d_i(x)= s_i(x)- {\underline s}(x),\endmatrix}\\
& h^{(\gamma) }(x)=\min_{A:\, \emptyset \ne A\subset \{1,2,\ldots,
n\}}\frac{1}{|A|} \bigg[\sum_{i\in A} \frac{d_i(x)}{q_i(x)^\gamma }+
\sum_{i\in A, j\notin A} \frac{|\partial_{ij} U (x)|}{[q_i(x)\vee
q_j(x)]^\gamma}
  \bigg], \\
  &\text{\hskip 8em} \gamma \ge 0,\tag 1.7\endalign$$
  where $a\vee b=\max\{a, b\}$ and $|A|$ is the cardinality of the set $A$.

\proclaim{Theorem 1.2} We have
$$\lambda_1 (U)\ge
\inf_{x\in \Bbb R^n} \bigg\{ {\underline
s}(x)+\frac{h^{(1/2)}(x)^2} {1+\sqrt{1- h^{(1)}(x)^2}} \bigg\}.
\tag 1.8 $$ \endproclaim

A close related topic to the Poincar\'e inequality is the logarithmic Sobolev inequality
with optimal constant $\sz(U)$:
$$\sz (U)\,\ent_{\mu}\big(f^2\big)\le 2 D(f), \qqd f\in {\scr D}(D), \tag 1.9 $$
where $\ent_{\mu}(f)=\mu(f\log f)-\mu(f)\log\mu(f)$ for $f\ge 0$.
Correspondingly, we have the conditional marginal inequality for
$\mu_U^{ x_{\backslash i}}$, given $x_{\backslash i}$,
 with optimal constant $\sz^{x_{\backslash i}}$:
$$ \sz^{x_{\backslash i}} \ent_{\mu_U^{ x_{\backslash i}}} (f)\le
\int_{\Bbb R} {f'}^2 \d \mu_U^{ x_{\backslash i}}, \qquad f\in
C_0^\infty (\Bbb R).   \tag 1.10$$

We can now state a very recent result due to [8; Theorem 1], which
is consistent with Theorem 1.1.

\proclaim{\thm\;1.3} The logarithmic Sobolev constant $\sz(U)\ge \lz_{\min}(A)$,
where the matrix $A=(A_{ij})$ is defined by
$$A_{ij}=\cases
\inf_x \sz^{x_{\backslash i}}, &\qd\text{if } j=i\\
-\sup_x |\partial_{ij}U(x)|. &\qd\text{if } j\ne i.
     \endcases$$
\endproclaim

In view of the above results, it is clear that the one-dimensional case plays a crucial role.
In that case, a representative result of the paper is as follows.

\proclaim {\prp\;1.4} In dimensional one, replace $U$ with $u_{\bz_1,\,\beta_2}^{}(x)=x^4-\beta_1 x^2+\beta_2 x$ for
some constants $\beta_1\ge 0$ and $\beta_2\in \Bbb R$. Then we have
$$\align
{4 e^{14}}\! \exp\bigg[\!-\frac 1 4 \bz_1^2+2\log(1+\bz_1)\bigg]&\ge \inf_{\bz_2}\lz_1(u_{\bz_1,\,\beta_2}^{})\\
&\ge \inf_{\bz_2}\sz(u_{\bz_1,\,\beta_2}^{})\\
&\ge  \frac{\sqrt{\beta_1^2+8}- \beta_1}{\sqrt{e}}
\exp\bigg[-\frac{1}{8}\beta_1\Big(
\beta_1+\sqrt{\beta_1^2+8}\Big)\bigg].
\endalign$$
In particular, $\inf_{\bz_2}\lz_1(u_{\bz_1,\,\beta_2}^{})$ and $\inf_{\bz_2}\sz(u_{\bz_1,\,\beta_2}^{})$
have the same order as \newline
$\exp[-\bz_1^2/4+O(\log \bz_1)]$ as $\bz_1\to\infty$.
\endproclaim

The exponent $\beta_1^2/4$ here equals, approximately as $\beta_1\to\infty$, the square of
the variance of a random variable having the distribution with density $\exp[-x^4+\beta_1 x^2]/Z$ on the real  line.

The remainder of the paper is organized as follows. In the next
section, we study an alternative variational formula for spectral
gap. This is especially meaningful in the context of diffusions. The
proofs of Theorems 1.1 and 1.2 are completed in \S 3. The
one-dimensional spectral gap is the main topic in \S 4. The
logarithmic Sobolev constant is studied in \S 5, in which \prp\;1.4
is proven. Even though the explicit and universal upper and lower
estimates, as well as the criteria, for the spectral gap and
logarithmic Sobolev constant are all known (cf. [13; Chapter 5,
Theorem 7.4] and \S 4 below), it is still quite a distance to arrive
at \prp\;1.4. Actually, we study this model several times (Examples
4.3, 4.6, 4.9, 5.3, and Proposition 4.7) by using different
approaches. Thus, a part of the paper is methodological, it takes
time and space to make some comparison of different methods. Two
typical infinite-dimensional models are treated in the last section.

\subhead{2. Alternative variational formula for spectral gap}
\endsubhead
Let $(E, {\scr E}, \mu)$ be a probability space and $ L^2(\mu)$ be
the ordinary $L^2$-space of real functions. Corresponding to a
$\mu$-reversible Markov process with transition probability $ P(t,
x, \cdot)$, we have a positive, strongly continuous, contractive
and self-adjoint semigroup $ \{P_t\}_{t\ge 0}$ on $ L^2(\mu)$ with
generator $ (L, {\scr D}(L))$. Throughout this section, $ (\cdot,
\cdot)$ and $\|\cdot\|$ denote, respectively, the inner product
and the norm in $ L^2(\mu)$. By elementary spectral theory, we
have
$$\frac{1}{t}(f-P_t f, f)\uparrow \text{ some } D(f, f)=:
D(f)\le \infty   \text{ as } t\downarrow 0. \tag 2.1$$ Set ${\scr
D}(D)=\{f\in L^2(\mu): D(f)<\infty \}$  and define $D(f,
g)=(D(f+g)-D(f-g))/4$ for $f, g\in {\scr D}(D)$. Then, $(D, {\scr
D}(D))$ is a Dirichlet form. Moreover,
$$D(f, g)=-(L f, g), \qquad f, g \in {\scr D}(L). \tag 2.2 $$
The formula in (2.4) below goes back to [14].

\proclaim{Theorem 2.1} The spectral gap $\lambda_1 (L)$ is
described by the largest constant $\kappa $ in the following
equivalent inequalities.
$$\align
 \kappa \var_{\mu} (f)&\le D(f), \;\;\qquad f\in {\scr D}(D), \tag 2.3\\
 \kappa D(f)&\le \|L f \|^2,  \qquad f\in {\scr D}(L). \tag
2.4\endalign$$       \endproclaim

\demo{Proof} Let $\{E_{\alpha})\}_{\alpha\ge 0}$ be the spectral
representation of $L$. Then  $L=-\int_0^\infty \alpha\d E_\alpha$.
The optimal constant $\kappa$ in (2.3) is known to be
$\lambda_1=\lambda_1 (L)$. Note that
$$\align
\|L f\|^2&=(L f, L f)\\
&=(f, L^2 f)\\
&=\bigg(f, \int_0^\infty \alpha^2 \d E_{\alpha} f\bigg)\\
& =
\int_0^\infty \alpha^2 \d (E_{\alpha} f, f )\\
&=\int_{\lambda_1}^\infty \alpha^2 \d (E_{\alpha} f, f )\ge
\lambda_1 \int_{\lambda_1}^\infty \alpha \d (E_{\alpha} f, f
)\\
&=\lambda_1
\int_{0}^\infty \alpha \d (E_{\alpha} f, f )\\
&= \lambda_1 (f, -L f)\\
&= \lambda_1 D(f).
\endalign$$
Because the only inequality here cannot be improved, the largest
constant $\kappa$ in (2.4) is also equal to
$\lambda_1$.\qed\enddemo

\proclaim{Remark 2.2} {\rm Actually, it is known and is also easy to check
that $(2.3)$ is equivalent to the correlation inequality
$$\lz_1(L)|\cov_{\mu}(f, g)|\le \big(D(f) D(g)\big)^{1/2},\qqd f, g\in {\scr D}(D),
\tag 2.5$$
where $\cov_{\mu}(f, g)=\mu(f g)-\mu(f)\mu(g)$ and $\mu(f)=\int f\d\mu$.}
See the comment below Proposition 3.2 for a proof.\endproclaim

Before moving further, let us mention that the above proof also
works for the principal eigenvalue. In this case, $L 1\ne 0$ and
$\mu$ can be infinite. Then the principal eigenvalue $\lambda_0$
can be described by the following equivalent inequalities.
$$\align
&\tilde \kappa \|f\|^2 \le D(f), \qquad f\in {\scr D}(D),\\
&\tilde \kappa D(f) \le \|L f\|^2, \qquad f\in {\scr D}(L). \tag
2.6
\endalign$$

The formula (2.4) is especially useful for diffusion on Riemannian
manifolds. Thus, the next result is meaningful for a more general
class of diffusion in $\Bbb R^n$ by using a suitable Riemannian structure.

\proclaim{Corollary 2.3} Let $L= \Delta - \langle \nabla U, \nabla
\rangle$ for some $U\in C^\infty  (\Bbb R^n)$ with $Z:= \int_{\Bbb
R^n} e^{-U}\d x$ $< \infty $ and set $\mu(\d x)= e^{-U} \d x/Z$. Then
$$\|L f\|^2=\int_{\Bbb R^n} \bigg[\sum_{i, j}(\partial_{ij} f)^2+
\langle\hes (U)\nabla f, \nabla f\rangle \bigg] \d \mu, \qquad
f\in C_0^\infty\big(\Bbb R^n\big), \tag 2.7
$$
where $\langle\cdot, \cdot\rangle$ stands the usual inner product
in ${\Bbb R}^n$. In particular, we have
$$\lambda_1 (U)\ge \inf_{x\in \Bbb R^n}\lambda_{\min}(\hes (U)(x)), \tag 2.8 $$
where $\lambda_{\min}(M)$ is the minimal eigenvalue of the matrix $M$.
\endproclaim

\demo{Proof} The proof of (2.7) is mainly a use of integration by
parts formula. Because $L f=\sum_i \big(\partial_{ii} f -
\partial_{i} U \partial_{i} f\big)$, we have
$$
\langle \nabla f, \nabla L f\rangle = \sum_j \partial_{j}
f\sum_{i} \partial_{j} \big(\partial_{ii} f -
\partial_{i} U \partial_{i} f\big)
=\sum_{i, j} \partial_{j} f \big(\partial_{iij} f -
\partial_{ij} f \partial_{i} U -\partial_{i} f \partial_{ij}
U\big).
$$
Next,
$$\align
\frac{1}{Z}\int_{\Bbb R^n} \sum_{j} \partial_{j} f \sum_i
\big(\partial_{iij} f -
\partial_{ij} f \partial_{i} U\big) e^{-U}
&=\frac{1}{Z}\int_{\Bbb R^n}  \sum_{j} \partial_{j} f \sum_i
\partial_i\Big(\partial_{ij} f e^{-U}\Big)\\
&=-\frac{1}{Z}\int_{\Bbb R^n}  \sum_{i, j} \big(\partial_{ij}
f\big)^2  e^{-U}\\
&=-\int_{\Bbb R^n}  \sum_{i, j} \big(\partial_{ij} f\big)^2 \d \mu,
\quad f\in C_0^\infty\big(\Bbb R^n\big).
\endalign$$
Noting that $\mu$ is a probability measure and the diffusion
coefficients are constants, the Dirichlet form is regular (cf. [12;
condition (4.13)] for instance). Actually, the martingale problem
for $L$ is well posed. Thus, $L C_0^\infty \big(\Bbb R^n\big)\subset
C_0^\infty \big(\Bbb R^n\big)\subset {\scr D} (L)$, and so
$$\|L f\|^2=\int_{\Bbb R^n}  L f \cdot L f \d \mu
= -\int_{\Bbb R^n} \langle\nabla L f, \nabla f\rangle\d\mu,\qquad
f\in C_0^\infty\big(\Bbb R^n\big).$$ Combining these facts
together, we get (2.7).

To prove the last assertion, applying Theorem 2.1 and (2.7), we
get
$$\align
\lambda_1(L)&=\inf_{f\in {\scr D}
(L),\,f\ne\,\text{const}}\frac{\|Lf\|^2}{D(f)}\\
& =\inf_{f\in
C_0^\infty(\Bbb R^n),\,f\ne\,\text{const}}\frac{\|Lf\|^2}{D(f)}\\
&\ge \inf_{f\in C_0^\infty(\Bbb
R^n),\,f\ne\,\text{const}}{\int_{\Bbb R^n} \langle\hes
(U)\nabla f, \nabla f\rangle \d \mu}\Big/{D(f)}\\
& \ge \inf_{x\in \Bbb
R^n}\lambda_{\min}(\hes (U)(x)). \qed \endalign$$
\enddemo

\proclaim{Remark 2.4} {\rm Actually, under the assumption of
Corollary $2.3$, the Bakry-Emery criterion (cf. [14] or [7;
Corollary 1.6]) implies a stronger conclusion:
$$\sz(U) \ge \inf_{x\in \Bbb R^n}\lambda_{\min}(\hes (U)(x)). \tag 2.9$$
}\endproclaim

A simple counterexample for which (2.8) and (2.9) are not effective is the following.
This example also shows that (1.4) is an improvement of (2.8).

\proclaim{Example 2.5} Consider the two-dimensional case. Let
$$U(x)=x_1^4+x_2^4-\bz \big(x_1^2+x_2^2\big) +2J x_1 x_2 $$ with
constants $\bz\ge 0$ and $J\in {\Bbb R}$. Then $\inf_{x\in \Bbb
R^2}\lambda_{\min}(\hes (U)(x))\le 0$ and $U$ is not convex at
infinity, but $\lz_1(U)>0$ in a region of $(\bz, J)\subset {\Bbb
R}\times {\Bbb R}_+$.
\endproclaim

\demo{Proof} First, we have
$$\hes (U)(x)=\left(\matrix 12 x_1^2 -2\bz & 2J\\
2J & 12 x_2^2 -2\bz\endmatrix\right).$$ Because for the matrix
$$A=\left(\matrix c_1 &  2J\\
2J& c_2\endmatrix\right),$$ we have
$\lz_{\min}(A)=2^{-1}\big(c_1+c_2-\sqrt{(c_1-c_2)^2+16 J^2}\,\big)$.
Hence
$$\lambda_{\min}(\hes (U)(x))=2\min_{x_1, x_2}\Big\{3\big(x_1^2+x_2^2\big)
 -\bz-\sqrt{9 \big(x_1^2-x_2^2\big)^2+J^2}\,\Big\}.$$
Setting $x_1=x_2=0$, we get
$$\inf_{x\in \Bbb R^2}\lambda_{\min}(\hes (U)(x))\le -2(\bz + |J|)\le 0.$$
Next, since
$$\lim_{|x_1|\to\infty}\Big(3 x_1^2 -\bz-\sqrt{9 x_1^4+J^2}\,\Big)
=\lim_{z\to 0}\frac{3 -\sqrt{9 +J^2 z^2}}{z}-\bz=-\bz,$$ we have
$$\varliminf_{|x|\to\infty}\lambda_{\min}(\hes (U)(x))\le
\varliminf_{x_2=0,\;|x_1|\to\infty}\lambda_{\min}(\hes (U)(x))
\le -2\bz\le 0.$$
This means that $U$ is not convex at
infinity. The last assertion of the example is the one of the main aims of
this paper and it is even true in the higher dimensions (cf. Theorem
6.3  below).\qed\enddemo

\subhead{3. Proofs of Theorems 1.1 and 1.2 and some remarks}
\endsubhead
As a preparation, we prove a result which is an improvement of (2.8)
and [7; Proposition 3.1]. We adopt the notation given in \S1.

\proclaim{Proposition 3.1} We have
$$\lambda_1 (U)\ge \inf_{x\in \Bbb R^n}
\lambda_{\min} \big({\widetilde \hes}(U)(x)\big). \tag 3.1 $$
\endproclaim

\demo{Proof} First, applying Theorem 2.1 and (2.7) to the $i$th
marginal, we have
$$\int_{\Bbb R} \big[(\partial_{ii} f)^2 + (\partial_{ii}U)
(\partial_{i} f)^2\big]\d \mu_U^{x_{\backslash i}} \ge
\lambda_1^{x_{\backslash i}} \int_{\Bbb R} (\partial_{i} f)^2 \d
\mu_U^{x_{\backslash i}}, \qquad f\in C_0^\infty\big(\Bbb
R^n\big). \tag 3.2
 $$
Next, denote by $\hes_0(U)$ the symmetric matrix obtained from the
Hessen matrix $\hes (U)$ replacing the diagonal elements with
zero. Then, by (3.2), we have
$$\align
&\int_{\Bbb R^n} \bigg[\sum_{i, j}(\partial_{ij} f)^2+
\langle\hes (U)\nabla f, \nabla f\rangle \bigg] \d \mu_U^{} \\
&\ge \sum_{i}\int_{\Bbb R^n} \bigg\{ \int_{\Bbb R}
\big[(\partial_{ii} f)^2+ (\partial_{ii} U) (\partial_i f)^2\big]
\d\mu_U^{x_{\backslash i}}\bigg\}\d \mu_U^{} -\sum_{i}\int_{\Bbb
R^n}\big[(\partial_{ii} U) (\partial_i f)^2\big]
\d \mu_U^{}\\
&\qquad +\int_{\Bbb R^n}
\langle\hes (U)\nabla f, \nabla f\rangle \d \mu_U^{}\\
&\ge \sum_{i}\int_{\Bbb R^n} \bigg\{\lambda_1^{x_{\backslash i}}
\int (\partial_i f)^2 \d\mu_U^{x_{\backslash i}}\bigg\}\d \mu_U^{}
  +\int_{\Bbb R^n}
\langle\hes_0 (U)\nabla f, \nabla f\rangle \bigg] \d \mu_U^{} \\
&= \sum_{i}\int_{\Bbb R^n} \lambda_1^{x_{\backslash i}}
(\partial_i f)^2 \d \mu_U^{}
  +\int_{\Bbb R^n}
\big\langle\hes_0 (U)\nabla f, \nabla f\big\rangle \d \mu_U^{}\\
&  =
\int_{\Bbb R^n}
\langle{\widetilde\hes} (U)\nabla f, \nabla f\rangle \d \mu_U^{} \\
&\ge \inf_{x\in \Bbb R^n} \lambda_{\min} \big({\widetilde\hes}
(U)(x) \big)\int_{\Bbb R^n} |\nabla f|^2 \d \mu_U^{}, \qquad f\in
C_0^\infty\big(\Bbb R^n\big). \tag 3.3
\endalign$$
Now, the required assertion follows from the proof of the last
assertion of Corollary 2.3. \qed\enddemo

From the proof of Proposition 3.1, it is clear that the only
argument where we may lose somewhat is the first inequality of
(3.3), since the terms $\sum_{i\ne  j}(\partial_{i j} f)^2$ are
ignored there. Hence the estimate (3.1) is mainly meaningful if the
interactions are not strong. The interacting potentials considered
in this paper are rather simple; for general interactions, one needs
some ``block estimates'' which are not touched here, instead of the
``single-site estimates'' studied in this paper.

The shorthand of (3.1) is that
the minimal eigenvalue $\lambda_{\min} \big({\widetilde\hes} (U)
\big)$ may not be computable in practice. For this, we need the second
variational procedure. To do so, let $s= \min_i
\big\{\lambda_1^{x_{\backslash i}}-\sum_{j: j\ne i} |\partial_{ij}U|\big\}$ and define
$$
q_{ij}=
\cases
 |\partial_{ij}U|, \quad &\text{if } i\ne j\\
s- \lambda_1^{x_{\backslash i}}, \quad &\text{if } i= j.
\endcases$$ Then, $Q:=(q_{ij})$, depending on $x$, is a symmetric
$Q$-matrix, not necessarily conservative \big(i.e., $\sum_{j}
q_{ij}\le 0$\big).

\demo{Proof of Theorem $1.1$} The first estimate in (1.4) follows
from Proposition 3.1. Next, by [15; Theorem 1.1], we have
$$ \lambda_{\min}(-Q) \ge \sup_{w>0} \min_i\big[-{Q w}/{w} \big](i),
\tag 3.4 $$ where $Q w(i)=\sum_{j} q_{ij} w_j$. We remark that the
sign of the equality in (3.4) holds once $Q$ is irreducible (cf.
[15; Proposition 4.1]). Noting that for every symmetric matrix
$B=(b_{ij})$ with nonnegative diagonals and any vector $w$, we have
$$\langle w, B w\rangle =\sum_i b_{ii}w_i^2 + 2\sum_{i\ne j}b_{ij}w_iw_j
\ge \sum_i b_{ii}w_i^2 -2\sum_{i\ne j}|b_{ij}w_iw_j|=\big\langle |w|, {\widetilde B} |w|\big\rangle, $$
where ${\widetilde B}=(\tilde b_{ij}): \tilde b_{ii}=b_{ii},\; \tilde b_{ij}=-|b_{ij}|$
for $i\ne j$ and $|w|=(|w_i|)$. Letting $w^*$ be a vector with $\langle w^*,  w^*\rangle=1$ such that $\lz_{\min}(B)=\langle w^*, B w^*\rangle$,
it follows that
$$\lz_{\min}(B)\ge \big\langle |w^*|, {\widetilde B} |w^*|\big\rangle\ge
\lz_{\min}\big({\widetilde B}\big)\langle |w|^*,  |w|^*\rangle=\lz_{\min}\big({\widetilde B}\big).$$
Based on this fact and as an application of (3.4), we get
$$\align
\lambda_{\min} \big({\widetilde\hes} (U)(x) \big) &\ge
\lambda_{\min} \big(\text{\rm diag}(s)-Q \big)\\
&= s+ \lambda_{\min} (-Q )\\
&\ge s+ \max_{w>0} \min_i
\bigg[-s +\lambda_1^{x_{\backslash i}}-\sum_{j: j\ne i} q_{ij} w_j/w_i\bigg]\\
&= \max_{w>0} \min_i \bigg[\lambda_1^{x_{\backslash i}}-\sum_{j: j\ne
i} q_{ij} w_j/w_i\bigg]. \tag 3.5\endalign$$ Combining this with
the first estimate in (1.4), we get the second one in (1.4), and
so complete the proof of Theorem 1.1. \qed\enddemo

\demo{Proof of Theorem $1.2$} In the above proof, replacing
$\lambda_1^{x_{\backslash i}}$, $s$ and $q_{ii}$ with
$\eta_i^{x_{\backslash i}}$, ${\underline s}(x)$ and $q_{ii}(x)$,
respectively, but keep $q_{ij}\,(i\ne j)$ to be the same, we
obtain
$$\lambda_{\min} \big({\widetilde\hes} (U)(x) \big)
\ge {\underline s}(x)+ \lambda_{\min}(-Q(x)).
 $$
By Proposition 3.1, it suffices to estimate $\lambda_{\min}(-Q(x))$.
Note that $\lambda_{\min}(-Q(x))$ is nothing but the principal
(Dirichlet) eigenvalue of $Q(x)$, often denoted by
$\lambda_0(Q(x))$. Because $Q(x)$ is symmetric, and so its
symmetrizing measure is just the uniform distribution on
$\{1,2,\ldots, n\}$. Now the conclusion of Theorem 1.2 follows from
[16; Theorem 1.1] plus some computations.
 \qed\enddemo

We conclude this section with some remarks.

Let $d_i=-q_{ii}-\sum_{j\ne i} q_{ij} $. By setting $w_i=$
constant in (3.4), it follows that $\lambda_{\min} (-Q)\ge \min_i
d_i$. The sign of the equality holds if $(d_i)$ is a constant.
Otherwise, this well-known simplest conclusion is usually rough.
For instance, take
$$ Q=\pmatrix -1 \quad & 1 \quad & 0\\
1 \quad & -2 \quad & 1\\
0 \quad & 1 \quad & -3
\endpmatrix. $$
Then $\lambda_{\min}(-Q)=2-\sqrt{3}>0 $ \big(the equality of (3.4) is attained at
the positive eigenvector $w=\big(2+\sqrt{3}\,, 1+\sqrt{3}\,, 1\big)$ but $\min_i d_i=0$. This
shows that the use of the variational formula (3.4) is necessary
to produce sharper lower bounds.

When $\partial_{ij} U\le 0$ for all $i\ne j$, then
${\widetilde\hes} (U)= \text{\rm diag}(s)-Q$, and so the sign of
the first equality in (3.5) holds. In this case, the estimate
(3.5) is quite sharp, since so is (3.4). However, for general
$\partial_{ij} U\,(i\ne j)$, the lower bound in (3.5) may be less effective
but we do not have a variational formula as (3.4) in such a
general situation.

For a given symmetric matrix $B=(b_{ij})$ $\big({\widetilde \hes
}(U), \text{ for instance}\big)$, the classical variational
formula, which is especially powerful for upper bounds, is as
follows.
$$\align
\lambda_{\min}(B)&=\inf\bigg\{\sum_{i,j} w_i b_{ij} w_j:\; \sum_i w_i^2=1\bigg\}\\
&=\inf\bigg\{\sum_{i}\bigg(b_{ii}+\sum_{j: j\ne i} b_{ij}\bigg)
w_i^2-\frac{1}{2} \sum_{i, j}b_{ij}(w_j-w_i)^2: \sum_{i}w_i^2=1
\bigg\}.\\
&\tag 3.6\endalign
$$
For a given symmetrizable $Q$-matrix $(q_{ij})$ with symmetric
probability measure $\mu$, set
$$D(f)=\frac{1}{2}\sum_{i,j} \mu_i
q_{ij}(f_j-f_i)^2+\sum_i \mu_i d_i f_i^2,$$ where
$d_i=-q_{ii}-\sum_{j\ne i} q_{ij} $ as defined before. Then, an
alternative formula of (3.6), in terms of the Donsker-Varadhan's
theory of large deviations, goes as follows.
$$\align
\lambda_{\min}(-Q)&=\inf_{f}\bigg\{D(f): \sum_i \mu_i f_i^2=1\bigg\}\\
&=\inf_{\alpha\ge 0}\bigg\{D\big(\sqrt{\d \alpha/ \d \mu}\big): \sum_i \alpha_i =1\bigg\}\\
&=\inf_{\alpha\ge 0}\bigg\{I( \alpha)+ \sum_i \alpha_i d_i: \sum_i \alpha_i =1\bigg\}\\
&=\inf_{\alpha\ge 0}\bigg\{-\inf_{u>0}\sum_{i, j}\alpha_i q_{ij}(u_j-u_i)/u_i+ \sum_i \alpha_i d_i: \sum_i \alpha_i =1\bigg\}\\
&=\inf_{\alpha\ge 0 }\bigg\{\frac{1}{2} \sum_{i,
j}\big(\sqrt{\alpha_i q_{ij}}- \sqrt{\alpha_j q_{ji}}\big)^2+
   \sum_{i} \alpha_i d_i:  \sum_i \alpha_i =1\bigg\},\\
   & \tag 3.7
\endalign$$
where $I$ is the $I$-functional in the theory of large deviations.
Refer to [17; Proof of Theorem 8.17] for more details. In other
words, the large deviation principle provides an alternative
description of the classical variational formula, but not (3.4), for
which one needs a variational formula for the Dirichlet forms (cf.
[15]).

Finally, we remark that the proof of \prp\;3.1 can be also used in
the study of other inequalities. The details are omitted here since
they are not used subsequently (cf. [7]). The next one is a partial
extension of (2.5).

\proclaim{Proposition 3.2} Under the assumption of Corollary $2.3$, we have for every invertible, nonnegative
and diagonal matrix $D$, the largest constant $\kz$:
$$\align
\kz \big|\cov_{\mu_U^{}}(f, g)\big|&\le \bigg(\int |D\nabla f|^2\d\mu_U^{}\int |D^{-1}\nabla f|^2\d\mu_U^{} \bigg)^{1/2},\\
&\qqd f, g\in C_0^\infty\big({\Bbb R}^n\big) \tag 3.8\endalign $$
 satisfies
 $$\kz \ge \inf_x \hat{\lz}_{\min}\big(D \,\widetilde\hes (U)\, D^{-1}(x)\big),\tag 3.9$$
where $\hat{\lz}_{\min}(M)=\max\{c: M\ge c\, \text{\rm Id}\}$.\endproclaim

Before moving further, let us make some remarks about the proof of Proposition 3.2. Note that
$$\int |D\nabla f|^2\d\mu=\int \langle D^2\nabla f, \nabla f\rangle\d\mu
$$
which is the Dirichlet form corresponding to the diffusion operator with diffusion coefficients $D^2$ and
potential $U$. Denote by $\lz_1\big(D^2, U\big)$ the spectral gap of the last operator, then we have
$$
\align
\big|\cov_{\mu}(f, g)\big|^2&\le \var_{\mu}(f)\var_{\mu}(g)\\
&\le \frac{1}{\lz_1\big(D^2, U\big)\, \lz_1\big(D^{-2}, U\big)}
\int \langle D^2\nabla f, \nabla f\rangle\d\mu\int \langle D^{-2}\nabla g, \nabla g\rangle\d\mu.\\
&\tag 3.10\endalign$$
Hence, we obtain a lower bound of the optimal constant in (3.8):
$$\kappa\ge \sqrt{\lz_1\big(D^2, U\big)\, \lz_1\big(D^{-2}, U\big)}.\tag 3.11$$
The proof is quite natural. Furthermore, by setting $D$ to be the identity matrix, we obtain
(2.5) with sharp constant. However, the estimate (3.11) is usually not sharp in the general case. Note that
the sign of the last equality in (3.11) holds if $f$ and $g$ are the correspondent
eigenfunctions with respect to the operators, but the sign of the first equality in (3.11) holds
iff $f$ and $g$ are proportional almost surely (due to the use of the Cauchy-Schwarz inequality).
This can happen only if $D$ is trivial: all the
diagonals of $D$ are equal.

A better way to study (3.8) is using the semigroup's approach. Write
$$\aligned
\cov_\mu (f, g)&=\int (f-\mu(f))g \d\mu\\
&=-\int \bigg(\int_0^\infty \frac {\d}{\d t} P_t f\d t\bigg)g\d\mu\\
&=-\int_0^\infty\bigg(\int g L P_t f \d\mu\bigg)\d t\\
&=\int_0^\infty\bigg(\int \langle \nabla P_t f, \nabla g
\rangle\d\mu\bigg)\d  t.
\endaligned$$
Now, as a good application of the Cauchy-Schwarz inequality, we get
$$\big|\cov_\mu (f, g)\big|\le \bigg[\int_0^\infty \bigg(\int |D\nabla P_t f|^2 \d\mu\bigg)^{1/2}\d  t\bigg]
\bigg(\int |D^{-1}\nabla g|^2 \d\mu\bigg)^{1/2}.$$ The problem is
now reduced to study the decay of $\int |D\nabla P_t f|^2 \d\mu$ in
$t$ (cf. [7]).

Similarly to Proposition 3.1, as checked by Feng Wang in 2002, we
have the following result which improves (2.9), but may be weaker
than Theorem 1.3.

\proclaim{Proposition 3.3} Under the assumption of Corollary $2.3$, we have
$$\sz(U)\ge \inf_x \lz_{\min} \big(\,{\overline \hes} (U)\big), \tag 3.12$$
where
$$\big(\,{\overline \hes} (U)\big)_{ij}=
\cases
\zz^{x_{\backslash i}},\qd &\text{if } j= i\\
(\hes (U))_{ij}\qd &\text{if } j\ne i,
\endcases$$
$\zz^{x_{\backslash i}}$ is the optimal constant in the inequality
$$
\zz^{x_{\backslash i}}\int f (\partial_{i} \log f)^2 \d
\mu_{U}^{x_{\backslash i}} \le \int f\ggz^i (\log f) \d
\mu_{U}^{x_{\backslash i}}, \qqd 0\le f\in C_0^\infty \big({\Bbb
R}^n\big),\tag 3.13 $$ and
$$\ggz^i (f)=(\partial_{ii}f)^2 +(\partial_{ii} U)(\partial_{i}f)^2.$$\endproclaim

\subhead{4. One-dimensional case. Explicit estimates}
\endsubhead
The operator now becomes $L=\d^2/\d x^2 - u'(x)\d/\d x$. Write
$b(x)=-u'(x)$. Then $b$ must have a real root. Otherwise, without
loss of generality, let $u'\ge \varepsilon
>0$. Then $-u$ is strictly decreasing, and so
$$\infty> Z:=\int_{\Bbb R} e^{-u}\ge \int_{-\infty}^0 e^{-u}
>  e^{-u(0)} \int_{-\infty}^0 1=\infty,$$
which is a contradiction.

Unless otherwise stated, throughout this section, we
consider the operator
$$ L=a(x)\frac{\d^2}{\d x^2} +b(x)\frac{\d}{\d x} .$$ Assume
that $a\in C(\Bbb R)$, $a>0$ and $Z=\int_{\Bbb R}e^{C(x)}/a(x)<
\infty $, where $C(x)=\int_0^x b/a$. Define $\mu(\d x)=(Z a
(x))^{-1}e^{C(x)} \d x $. Recall that
$$\lambda_1(L)=\inf\{D(f):
f\in C^1(\Bbb R),\; \mu(f)=0,\; \mu(f^2)=1 \},$$ where $D(f)=\int_{\Bbb
R} a {f'}^2 \d \mu$.

Let $\theta $ be a fixed real root of $b$. Choose $K=K_{\theta
}\in C(\Bbb R\setminus \{\theta \})$ such that $K$ is increasing (i.e., non-decreasing)
in $x$ when $|x-\theta |$ increases, $K(\theta \pm 0)>-\infty $,
and moreover
$$ K(r)\le \inf_{x:\, \pm(x-r)> 0}\big[-b(x)/(x-\theta )\big]
\qquad \text{for all $\pm(r-\theta )> 0$}, \tag 4.1$$ where and in
what follows, the notation ``$\pm $'' means that there are two
cases: one takes ``$+$'' (resp.,  ``$-$'') everywhere in the statement. Define
$$ \align
&\text{\hskip-2em}F(s)=F^r(s)=\int_{\theta }^s \frac{u-\theta
}{a(u)}\big[K(r)-K(u)\big] \d u,
\qquad  s,\; r\in \Bbb R, \tag 4.2\\
&\text{\hskip-2em}\delta_{\pm}(K)= \sup_{r:\, \pm(r-\theta )>0} K(r)
\inf_{s: \, \pm(r-\theta )>\pm (s-\theta )>0}
\frac{(s-\theta )\exp[-F(s)]}{\int_{\theta}^s \exp[-F(u)]\d u}\tag 4.3\\
&\text{\hskip-2em}\qquad\quad  \ge  \sup_{r:\, \pm(r-\theta )>0}
K(r) \exp[-F(r)]. \tag 4.4
\endalign$$

The next result is a modification of [12; Corollary  3.5]. It is
specially useful for those $b$ growing at least linear.

\proclaim{Theorem 4.1} \roster
\item By using the above notations, we have
$$\lambda_1(L) \ge \delta_+(K) \wedge \delta_-(K). \tag 4.5  $$
\item Suppose additionally that $K$ is a piecewise $C^1$-function, then we have
$$ \delta_{\pm} (K)\ge
K(r_{\pm}^{}) \exp\bigg[-\int_{\theta }^{r_{\pm}^{}}
\bigg(\int_{\theta }^x \frac{u-\theta }{a(u)}\d u \bigg)\d K(x)
\bigg], \tag 4.6$$ where
$$ r_{\pm}^{} =\pm \infty, \;\;\qquad \;\, \text{if } \lim_{r\to \pm
\infty } K(r) \int_{\theta }^{r } \frac{u-\theta }{a(u)}\d u \le
1, \tag 4.7
$$
and otherwise, $r_{\pm}^{}$ is the unique solution to the equation
$$K(r) \int_{\theta }^r
\frac{u-\theta }{a(u)}\d u =1, \qquad  \pm(r-\theta )>0. \tag
4.8$$
\endroster
\endproclaim

\demo{Proof} (a) First, consider the half-line $(\theta , \infty
)$. Assume that $K(r_{1}^{})>0$ for some $r_{1}^{}\in (\theta ,
\infty )$. Otherwise, (4.5) becomes trivial. Fix $r=r_{1}^{}$ and
define
$$ f_+(x)=\int_{\theta }^x \d y \exp[-F(y\wedge r_{1}^{})],\qquad x\ge \theta. $$
Then, we have $f_+>0$ on $(\theta, \infty),\; f_+(\theta)=0$, $ f_+'(\theta)=1$ and
$$\align
&f_+'(x)= \exp[-F(x\wedge r_{1}^{})]>0,\\
&f_+''(x)= - \frac{x-\theta }{a(x)} \big[ K(r_{1}^{}) -K(x\wedge
r_{1}^{})\big] f_+'(x)\le 0, \qquad x\ge \theta .
\endalign$$
Since $a\in C(\Bbb R)$, $a>0$, $K\in C(\Bbb R\setminus \{\theta
\})$ and $K(\theta +)$ is finite, we have $f_+\in C^2(\theta,
\infty)$.

Next, because $K$ is increasing on $(\theta , \infty )$
and $K(x)\le -b(x)/(x-\theta )$ for all $x>\theta $, we have
$$\align
-(a f_+''+ b f_+')(x)
&=\big\{(x-\theta )[K(r_{1}^{})-K(x\wedge r_{1}^{})]- b(x)\big\}f_+'(x)\\
&\ge \big\{(x-\theta )K(r_{1}^{})- (x-\theta )K(x)- b(x)\big\}f_+'(x)\\
&\ge (x-\theta )K(r_{1}^{}) f_+'(x), \qquad x>\theta . \tag 4.9
\endalign$$ Since $f_+''\le 0$, $f_+'$ is decreasing. By the Cauchy mean
value theorem, it follows that $(x-\theta )/f_+(x)$ is increasing
on $(\theta , \infty )$. Hence, by (4.9), we obtain
$$\align
-\bigg[\frac{ a f_+''+ b f_+'}{ f_+}\bigg](x)
&\ge \frac{r_{1}^{}-\theta }{f_+(r_{1}^{})} K(r_{1}^{}) f_+'(x) =
\frac{r_{1}^{}-\theta }{f_+(r_{1}^{})} K(r_{1}^{})
f_+'(r_{1}^{}),\\
&\qquad x\ge r_{1}^{}. \tag 4.10\endalign
$$ Combining (4.9) with (4.10), it follows that
$$\inf_{x>\theta }\bigg[ -\frac{ a f_+''+ b f_+'}{ f_+} \bigg] \ge K(r_{1}^{}) \inf_{s\in (\theta , r_{1}^{})}
\frac{(s-\theta )f_+'(s)}{f_+(s)}. $$
By (4.3), we have thus obtained
$$ \inf_{x>\theta } \bigg[-\frac{ a f_+''+ b f_+'}{ f_+}\bigg]\ge \delta_+(K).
\tag 4.11$$

(b) Next, consider the half-line $(-\infty , \theta )$. The proof
is parallel to (a). Let $K(r_{1}^{})>0$ for some $r_{1}^{}<\theta
$. Fix $r=r_{1}^{}$ and define
$$ f_-(x)=\int_{\theta }^x \d y \exp[-F(y\vee r_{1}^{})],\qquad x\le \theta. $$
Then $f_-<0$ on $(-\infty, \theta)$, $f_-(\theta )=0$, $f_-'>0$,
$f_-'(\theta )=1$ and $$f_-''(x)=-\frac{x-\theta
}{a(x)}[K(r_{1}^{})-K(x\vee r_{1}^{})]f_-'(x)\ge 0 $$ for all $x\le
\theta $. Moreover $f_-\in C^2(-\infty, \theta )$. Then
$$\align
-(a f_-''+ b f_-')(x)
&=\big\{(x-\theta )[K(r_{1}^{})-K(x\vee r_{1}^{})]- b(x)\big\}f_-'(x)\\
&\le \big\{(x-\theta )K(r_{1}^{})- (x-\theta )K(x)- b(x)\big\}f_+'(x)\\
&\le (x-\theta )K(r_{1}^{}) f_-'(x), \qquad x < \theta.
\endalign$$ Since $f_-<0$ and $f_-''\ge 0$, we have
$$-\bigg[\frac{ a f_-''+ b f_-'}{ f_-}\bigg](x)
\ge \frac{r_{1}^{}-\theta }{f_-(r_{1}^{})} K(r_{1}^{}) f_-'(x) =
\frac{r_{1}^{}-\theta }{f_-(r_{1}^{})} K(r_{1}^{})
f_-'(r_{1}^{}),\qquad x\le r_{1}^{}.
$$ Combining the last two inequalities with (4.3), we get
$$\inf_{x <\theta }\bigg[ -\frac{ a f_-''+ b f_-'}{ f_-}\bigg]
\ge \sup_{r_{1}^{}<\theta } K(r_{1}^{}) \inf_{s\in (\theta ,
r_{1}^{})} \frac{(s-\theta )f_-'(s)}{f_-(s)}= \delta_-(K).
$$

Finally, let $f=f_+ I_{[\theta , \infty )}+ f_- I_{(- \infty,
\theta )}$. Then $f\in C^2(\Bbb R)$ and
 $$\align
 \inf_{x \ne \theta } \bigg[-\frac{ a f''+ b f'}{ f}\bigg](x)
&=\bigg[\inf_{x >\theta } -\frac{ a f_+''+ b f_+'}{ f_+}(x)
\bigg]\wedge
\bigg[\inf_{x <\theta } -\frac{ a f_-''+ b f_-'}{ f_-}(x)
\bigg]\\
&\ge \delta_+(K)\wedge \delta_-(K).\endalign
 $$
The estimate (4.5) now follows from the last assertion of [12;
Theorem 3.1].

(c) To prove (4.4), noticing that $K$ is monotone, we may apply
the integration by parts formula and rewrite $F$ as follows.
$$\align
F(r) &=\int_{\theta }^r \frac{u-\theta}{a(u)}[K(r)-K(u)]\d u  \\
&=K(r) \int_{\theta }^r \frac{u-\theta}{a(u)} \d u
-\int_{\theta }^r K(u) \d\bigg(\int_{\theta }^u \frac{z-\theta}{a(z)}\d z\bigg)  \\
&=\int_{\theta }^r K'(u) \bigg(\int_{\theta }^u \frac{z-\theta}{a(z)} \d
z\bigg)\d u, \qquad r\ne \theta . \tag 4.12 \endalign$$ By the
assumption on $K$, it follows that $F\ge 0$, $F(r)$ is increasing in
$r$ as $|r-\theta |$ increases. Hence, by (4.2), we have
$$\align
\delta_{\pm}(K)&\ge \sup_{r:\, \pm(r-\theta )>0} K(r) \inf_{s:\,
\pm(r-\theta )> \pm (s-\theta )>0}\exp[-F(s)]\\
& =\sup_{r:\,
\pm(r-\theta )>0} K(r) \exp[-F(r)].\endalign $$ The proof of (4.4) is done.

(d) The second part of the theorem is to compute $\sup_{r\ne \theta} G(r)$, where $G(r)= K(r)
\exp[-F(r)]$. The answer is given by (4.6). To do so, first consider
the half line $(\theta, \infty)$.
Because $K$ is a piecewisely $C^1$, we may assume that $\big(\theta, \infty\big)=\cup_i (c_i, d_i]$,
$K\in C^1(c_i, d_i)$ and $K'\ge 0$ on $(c_i, d_i)$ for every $i$. By (4.12), we have for every $i$,
$$G'(r)=K'(r)\bigg[1- K(r)\int_{\theta }^r
\frac{u-\theta }{a(u)}\d u \bigg] \exp[-F(r)],\qquad r\in (c_i, d_i). \tag 4.13$$
Let
$\lim_{r\to\infty}K(r) \int_{\theta }^r \frac{u-\theta }{a(u)}\d u
\le 1$ and set $\tilde\theta_+=\inf\big\{r> \theta: K(r)>0\big\}$.
Note that $K$ is increasing, $K>0$ on $\big(\tilde\theta_+, \infty\big)$, and so
$K(r) \int_{\theta }^r \frac{u-\theta }{a(u)}\d
u $ is strictly increasing on $\big(\tilde\theta_+, \infty\big)$,
but is less or equal to zero on
$\big(\theta, \tilde\theta_+\big)$ when $\theta< \tilde\theta_+$.
It follows that $K(r) \int_{\theta }^r \frac{u-\theta }{a(u)}\d
u \le 1$ for all $r\in (\theta, \infty)$.
By (4.13), we have $G'(r)\ge 0$ on
every $(c_i , d_i)$ since so does $K'(r)$. This fact plus the continuity of $G$
implies that
$\sup_{r>\theta}G(r)=\lim_{r\to \infty }G(r)$.

Otherwise, we have
$$\align
& \lim_{r\to\infty}K(r) \int_{\theta }^r \frac{u-\theta }{a(u)}\d u
> 1 \quad \text{and }\\
&{\lim_{r\to{\tilde\theta_+}+}K(r)\cases
=0<\Big( \int_{\theta}^{\tilde\theta_+} \frac{u-\theta }{a(u)}\d
u \Big)^{-1}\quad \text{if } \theta<\tilde\theta_+\\
< \infty=\Big( \int_{\theta}^{\theta} \frac{u-\theta }{a(u)}\d
u \Big)^{-1}
\quad \text{if } \tilde\theta_+=\theta.
\endcases}
\endalign$$
Since $K(r)$ is increasing and $ \Big(\int_{\theta }^r \frac{u-\theta }{a(u)}\d
u \Big)^{-1} $\! is strictly decreasing,  the
curves $K(r)$ and $ \Big(\int_{\theta }^r \frac{u-\theta }{a(u)}\d
u \Big)^{-1} $ must have uniquely an intersection on $\big(\tilde\theta_+, \infty
\big)$, or equivalently on $(\theta , \infty
)$. So we have $\sup_{r>\theta} G(r)=\sup_{r>\tilde\theta_+ } G(r)= G(r_{+}^{})$, where
$r_{+}^{}$ is the unique solution to the equation (4.8).

The proof of the assertions on $(-\infty, \theta)$ is parallel.
 \qed\enddemo

The next two examples illustrate the applications of Theorem 4.1, and are treated several times in the paper.

\proclaim{Example 4.2} Let $u(x)=\alpha x^2+\beta x$ for some
constants $\alpha >0$ and $\beta\in \Bbb R$, and $a(x)\equiv 1$. Then we
have $\lambda_1(L_{\alpha,\,\beta})\ge \delta_+(K)\wedge
\delta_-(K)= 2\alpha $ which is exact.
\endproclaim

\demo{Proof} Since $-b(x)= -2\alpha x-\beta  $, we have root
$\theta =-\beta /(2\alpha )$, and so ${-b(x)}/{(x-\theta
)}$ $=2\alpha $. Thus, $K(r)=$ constant $2\alpha $. By (4.3), we get
$\delta_{\pm}(K)=2\alpha $ as claimed. It is easy to check that
the estimate is exact, since the corresponding eigenfunction is
linear.
 \qed\enddemo

\proclaim{Example 4.3} Let $u(x)=x^4-\beta_1 x^2+\beta_2 x$ for
some constants $\beta_1, \beta_2\in \Bbb R$ and $a(x)\equiv 1$. Then we
have
$$
\lambda_1(L_{\beta_1,\,\beta_2})\ge \delta_+(K)\wedge
\delta_-(K)\ge \frac{ \sqrt{\beta_1^2+2}- \beta_1}{\sqrt{e}} \exp\bigg[
-\frac{1}{2}\beta_1 \Big(\beta_1+ \sqrt{\beta_1^2+2}\,\Big)\bigg]
$$
uniformly in $\beta_2$. When $\beta_2=0$, we have
$$\lambda_1(L_{\beta_1,\,\beta_2})\ge \delta_+(K)\wedge
\delta_-(K)\ge \frac{ \sqrt{\beta_1^2+8}- \beta_1}{\sqrt{e}} \exp\bigg[
-\frac{1}{8}\beta_1 \Big(\beta_1+ \sqrt{\beta_1^2+8}\,\Big)\bigg].$$
\endproclaim

\demo{Proof} First, we have $b(x)=-u'(x)=-4 x^3 + 2\beta_1
x-\beta_2$. Let $\theta $ be a  real root of $u'$. For instance, we may take
$$\theta=\cases
0 &\quad \text{if } \beta_2=0\\
-\Big(\frac{\beta_2}{4}\Big)^{1/3} &\quad \text{if } \beta_1=0\\
2\sqrt{\frac{-\beta_1}{6}}\,\sinh\!\Big(\frac 1 3 \text{arc\,sinh}\, C\Big) &\quad \text{if } \beta_1<0\\
2\sqrt{\frac{\beta_1}{6}}\,\text{sgn}\,(C)
\cosh\!\Big(\frac 1 3 \text{arc\,cosh}\big( \text{sgn}\,(C)\,C\big)\Big) &\quad \text{if } \beta_1>0\text{ and }|C|>1\\
2\sqrt{\frac{\beta_1}{6}}\,
\cos\!\Big(\frac{4}{3}\pi+\frac 1 3 \text{arc\,cos}\, C\Big) &\quad \text{if } \beta_1>0\text{ and }|C|\le 1,
\endcases$$
where $C=\beta_2\Big(\frac{3}{2|\beta_1|}\Big)^{3/2}$. The reason we choose $4\pi/3$ rather
than $0$ or $2\pi/3$ in the last line is for the consistency of the case
 $\beta_2=0$. However, in what follows, we will not use the explicit formula of $\uz$, we are going to
 work out only the estimate uniform in $\uz$. Because
$$\frac{-b(x)}{x-\theta } = 4(x-\theta)^2 + 12 \theta (x-\theta)+12 \theta^2
-2\beta_1=
4(x+\theta/2)^2 +3 \theta^2 -2\beta_1,$$ we obtain
$$\align
&{\inf_{x>r}\frac{-b(x)}{x-\theta }= \cases
4(r+\theta /2)^2 + 3\theta^2-2\beta_1,\quad & \text{if }\; r\ge -\theta /2 \\
3\theta^2 - 2\beta_1, \quad &  \text{if }\; r\le -\theta /2,
\qquad r\ge \theta.  \endcases}\\
&{\inf_{x<r}\frac{-b(x)}{x-\theta }= \cases
4(r+\theta /2)^2 + 3\theta^2-2\beta_1,\quad & \text{if }\; r\le -\theta /2 \\
3\theta^2 - 2\beta_1, \quad &  \text{if }\; r\ge -\theta /2,
\qquad r< \theta.  \endcases}\endalign
$$
Naturally, one may define $K(r)$ as the right-hand sides, but then
the computations for the lower bounds of $\delta_{\pm}(K)$ become
very complicated. Here, we adopt a simplification.  Set $r_{\theta
}^{}=r-\theta $. Because
$$\align
&4(r+\theta /2)^2 + 3\theta^2-2\beta_1=12
(\theta +r_{\theta }^{}/2)^2+r_{\theta }^2- 2\beta_1 \ge r_{\theta
}^2 - 2\beta_1,\\
& 3\theta^2-2\beta_1\ge 9\theta^2/4-2\beta_1,\endalign$$
when $r\ge \theta $ (equivalently, $r_{\theta }^{}\ge 0$), we can
choose
$$
K(r)=K_{\theta}(r)=\cases
r_{\theta }^2-2\beta_1,\quad &\text{if }\; r_{\theta}^{} > -3\theta /2 \\
9{\theta }^2/4-2\beta_1,\quad &\text{if }\; r_{\theta }^{}<
-3\theta /2.  \endcases$$ By symmetry, one can define $K(r)$ for
the case of $r\le \theta $ as follows:
$$
K(r)=\cases
r_{\theta }^2-2\beta_1,\quad &\text{if }\; r_{\theta}^{} < -3\theta /2 \\
9{\theta }^2/4-2\beta_1,\quad &\text{if }\; r_{\theta }^{} >
-3\theta /2.  \endcases$$ Obviously, $K$ is a continuous piecewise $C^1$-function.

Suppose that $\theta <0$ for a moment. We use the notation $G(r)$
defined in the proof (d) of Theorem 4.1. Since $G(r)$ is
continuous in $r$, $G(r)$ is equal to the constant $K(-\theta/2)$
on $(\theta, -\theta/2]$, and $K'>0$ on $(\theta, \infty)$,
we have $\sup_{r>\theta} G(r)= \sup_{r\ge -\theta/2} G(r)
$. Clearly, $\lim_{r\to\infty}K(r)\int_{\theta}^r (u-\theta)\d u=\infty$ and hence we
can ignore (4.7) and handle with (4.8) only. There are two cases.

(a) Let $K(-\theta /2)\int_{\theta }^{-\theta /2} (u-\theta )\d
u<1$. That is $9 \theta^2/4<  \beta_1+\sqrt{
\beta_1^2+2}$. In this case, the solution to (4.8) should satisfy $r_+-{\theta} > -3\theta/2$. Solving equation
$$\big( r_{\theta}^2 -2\beta_1\big)\int_{\theta}^{r} (u-\theta)\d u = 1 ,
\qquad r_{\theta}^{}>-3\theta/2,$$ we get $(r_{+}^{}-\theta)^2=\beta_1+
\sqrt{\beta_1^2+2}$. Then
$$\align
-\frac{1}{2}\int_{\theta }^{r_{+}^{}} (x-\theta )^2 K'(x) \d x
& =
-\frac{1}{2}\int_{-3\theta/2}^{r_{+}^{}-\theta } x^2 \cdot 2 x \d
x \\
&= -\frac{1}{4} (r_{+}^{} -\theta)^4 + \frac{81}{64} \theta^4\\
&=-\frac{1}{4}\Big( \sqrt{\beta_1^2+2} + \beta_1\Big)^2 +
\frac{81}{64} \theta^4.\endalign $$ Hence we obtain
$$\align
 \sup_{r\ge -\theta/2} G(r)
 &= G(r_{+}^{})
 \ge \Big( \sqrt{\beta_1^2+2}- \beta_1\Big)
\exp\bigg[-\frac{1}{4}\Big( \sqrt{\beta_1^2+2}+ \beta_1\Big)^2 +
\frac{81}{64} \theta^4\bigg]\\
&\ge \Big( \sqrt{\beta_1^2+2}- \beta_1\Big)
\exp\bigg[-\frac{1}{4}\Big( \sqrt{\beta_1^2+2}+ \beta_1\Big)^2
\bigg]\\
&=\frac{ \sqrt{\beta_1^2+2}- \beta_1}{\sqrt{e}} \exp\bigg[
-\frac{1}{2}\beta_1 \Big(\beta_1+ \sqrt{\beta_1^2+2}\,\Big)\bigg]. \endalign$$

(b) Let $K(-\theta /2)\int_{\theta }^{-\theta /2} (u-\theta )\d
u\ge 1$. Equivalently, $9 \theta^2/4\ge  \beta_1+\sqrt{
\beta_1^2+2}$. In this case, the solution to (4.8) satisfies $r_+\in (\theta, -\theta/2)$.
Since $K$ is a constant on $(\theta, -\theta/2)$, by (4.13) and (4.12), $G=K$ on $(\theta, -\theta/2]$.
Hence
$$\align
\sup_{r> \theta} G(r) &=G(r_+)=K(-\theta/2)=\frac{9}{4}\theta^2 -
2\beta_1
\ge \sqrt{\beta_1^2+2}\,-\beta_1\\
&\ge  \Big( \sqrt{\beta_1^2+2}\,- \beta_1\Big)
\exp\bigg[-\frac{1}{4}\Big( \sqrt{\beta_1^2+2}\,+ \beta_1\Big)^2
\bigg]. \endalign$$

Combining (a) with (b) and (4.6), we obtain
$$\delta_+(K)\ge \frac{ \sqrt{\beta_1^2+2}- \beta_1}{\sqrt{e}} \exp\bigg[
-\frac{1}{2}\beta_1 \Big(\beta_1+ \sqrt{\beta_1^2+2}\,\Big)\bigg]. $$

Next, we estimate $\delta_-(K)$. Now, $K(r)=r_{\theta }^2- 2\beta_1$ on $(-\infty , \theta )$
since $\theta<0$.
From (4.8), we get the same solution
$(r_- -\theta)^2=\beta_1+ \sqrt{\beta_1^2+2}$. But
$$\align
 -\frac{1}{2}\int_{\theta }^{r_-} (x-\theta )^2 K'(x) \d x &=
-\frac{1}{2}\int_{0}^{r_- -\theta } x^2 \cdot 2 x \d x \\
&=
-\frac{1}{4} (r_- -\theta)^4 =-\frac{1}{4}\Big( \sqrt{\beta_1^2+2}
+ \beta_1\Big)^2 . \endalign$$ By (4.6) again, we get $$ \delta_-(K)\ge
\frac{ \sqrt{\beta_1^2+2}\,- \beta_1}{\sqrt{e}} \exp\bigg[
-\frac{1}{2}\beta_1 \Big(\beta_1+ \sqrt{\beta_1^2+2}\,\Big)\bigg].$$ Therefore, we have
proved the required lower bound in the case of $\theta <0$.

By symmetry, the same conclusion holds when $\theta >0$. The proof
for $\theta=0$ is much simpler as shown below.

When $\beta_2=0$, we simply let $\theta=0$. Then
$$\frac{-b(x)}{x}=4 x^2 - 2 \beta_1, \qquad x\ne 0.$$
We choose $K(r)= 4 r^2 - 2 \beta_1$. Then the equation (4.8) gives us
$$r_{\pm}^2=\frac 1 4 \Big(\beta_1 + \sqrt{\beta_1^2+8}\,\Big).$$
Because
$$\int_0^r \bigg[\int_0^x u \d u\bigg]\d K(x)= \int_0^r4 x^3 \d x=r^4,$$
by (4.6), we obtain the last required assertion.
\qed\enddemo

We will improve the estimate of Example 4.3 in \S 5 (Example 5.3) by a different method.

Before moving further, let us make some remarks about the estimate given
in Example 4.3. Recall that at the beginning of the proof,
in choosing the function $K(r)$, the term $12(\theta +r_{\theta}/2)^2$ was
removed, this simplified greatly the proof since the original quartic
equation is reduced to a quadratic one. For this reason one may worry
lost too much in the estimation and we want to know the best estimate we can get
by part (2) of Theorem 4.1. For this, we use a different trick. Consider the case of $\theta<0$ only.
We use the complete form of $K$:
$$K(r)=
\cases 4 (r+\theta/2)^2 + 3\theta^2 - 2 \beta_1, &\quad\text{if } r\ge -\theta/2\\
3\theta^2 -2\beta_1,  &\quad\text{if } \theta\le r\le -\theta/2\\
4 (r+\theta/2)^2 + 3\theta^2 - 2 \beta_1, &\quad\text{if } r< \theta.
\endcases$$

(i) Following the proof of Example 4.3, we study first the estimation of
$\delta_+(K)$. There are two cases.

(a) Let $K(-\theta/2) \int_{\theta}^{-\theta/2}(u-\theta) \d u<1$. That is, $3\theta^2< \beta_1 + \sqrt{\beta_1^2+8/3}$.
The idea is that in looking for a uniform estimate, we may
regard $r$ as a parameter rather than $\theta$. In other words, instead of solving equation
(4.8)
$$\big(4 r_{\theta}^2 + 12 \theta r_{\theta} + 12\theta^2-2 \beta_1\big) \int_{\theta}^r(u-\theta)\d u=1,
\qquad r_{\theta}>-3\theta/2$$
in $r$, we solve the equation in $\theta$. Then the equation
has two solutions:
$$\theta=\frac 1 6 \Big(-3r_{\theta} \pm \sqrt{6\beta_1 + 6/r_{\theta}^2 -3 r_{\theta}^2}\,\Big).$$
Since $\theta$ is real, $r_{\theta}$ must satisfy
$$r_{\theta}^2\le \beta_1 + \sqrt{\beta_1^2+2}. \tag 4.14$$
Next, in the ``$+$'' case, $\theta <0$ iff
$$r_{\theta}^2> \Big(\beta_1 + \sqrt{\beta_1^2+8}\,\Big)\Big/{4}, \tag 4.15$$
and it is obvious that $r_{\theta}>-3\theta/2$.
In the ``$-$'' case, it is automatically that $\theta<0$ and $r_{\theta}>\-3\theta/2$ iff
$$r_{\theta}^2> \Big(3\beta_1 + \sqrt{9\beta_1^2+24}\,\Big)\Big/{4}, \tag 4.16$$
To estimate the decay exponent, note that on the one hand, we have
$$\align
\theta r_{\theta}^3&=\frac 1 6 r_{\theta}^2 \Big(-3r_{\theta}^2 \pm \sqrt{6\beta_1 r_{\theta}^2 + 6 -3 r_{\theta}^4}\,\Big)\\
&=\frac 1 6 z \Big(-3z \pm \sqrt{6\beta_1 z + 6 -3 z^2}\,\Big),
\endalign$$
where $z=r_{\theta}^2$. On the other hand, we have
$$-\frac 1 2 \int_{-3\theta/2}^{r_+-\theta} x^2 K'(x)\d x
=- r_{\theta}^4 - 2\theta r_{\theta}^3 -\frac{27}{16}\theta^4.$$
Replacing $r_{\theta}^2$ with $z$ on the right-hand side plus
some computation, we finally get
$$\align
-\frac 1 2 \int_{-3\theta/2}^{r_+-\theta} x^2 K'(x)\d x
&= -\frac{3}{64}\bigg[-2 z^2 +8\beta_1 {z}+\beta_1^2+8 +\frac{2\beta_1}{z}+\frac 1 {z^2}\bigg]\\
&\quad \pm \frac{\sqrt{3}}{96}\,\big(9+9 \beta_1 z -23 z^2\big)
\sqrt{\frac{2\beta_1}{z}+\frac{2}{z^2}-1}.
\endalign$$
To obtain the uniform lower bound, by (4.14) and (4.15), we need to minimize the right-hand side under
the constrain
$$\cases
\big(\beta_1 +\sqrt{\beta_1^2+8}\,\big)/4 < z\le \beta_1 +\sqrt{\beta_1^2+2} \quad & \text{in the ``$+$'' case}\\
\big(3\beta_1 + \sqrt{9\beta_1^2+24}\,\big)\big/{4}< z\le \beta_1 +\sqrt{\beta_1^2+2} \quad & \text{in the ``$-$'' case}.
\endcases$$
A numerical computation shows that the first case is smaller than the second one and its
leading term is approximately $- 0.8\, \beta_1^2$.

(b) Let $K(-\theta/2) \int_{\theta}^{-\theta/2}(u-\theta) \d u\ge 1$. That is, $3\theta^2\ge \beta_1 + \sqrt{\beta_1^2+8/3}$.
Then we have the lower bound $\sqrt{\beta_1^2+8/3}-\beta_1$ which is decayed slowly than exponential.

(ii) Next, in the case of $r<\theta$, by assumption, $\theta<0$ and $r_{\theta}<0$, we have only
one solution
$$\theta=\frac 1 6 \Big(-3r_{\theta} - \sqrt{6\beta_1 + 6/r_{\theta}^2 -3 r_{\theta}^2}\,\Big), $$
and furthermore $\theta<0$ iff
$$r_{\theta}^2 < \Big(\beta_1 + \sqrt{\beta_1^2+8}\,\Big)\Big/{4}. $$
To estimate the decay exponent, note that on the one hand, since $r_{\theta}<0$, we have
$$\align
\theta r_{\theta}^3&=\frac 1 6 r_{\theta}^2 \Big(-3r_{\theta}^2 + \sqrt{6\beta_1 r_{\theta}^2 + 6 -3 r_{\theta}^4}\,\Big)\\
&=\frac 1 6 z \Big(-3z + \sqrt{6\beta_1 z + 6 -3 z^2}\,\Big).
\endalign$$
On the other hand, we have
$$-\frac 1 2 \int_{0}^{r_--\theta} x^2 K'(x)\d x
=- r_{\theta}^4 - 2\theta r_{\theta}^3.$$
Hence
$$-\frac 1 2 \int_{0}^{r_--\theta} x^2 K'(x)\d x
= -\frac{1}{\sqrt{3}}\, z \sqrt{2\beta_1 z + 2 - z^2}.$$
To obtain the uniform lower bound, it suffices to minimize the right-hand side under
the constrain
$$0< z<
\Big(\beta_1 +\sqrt{\beta_1^2+8}\,\Big)\Big/4 .$$
A numerical computation shows that the resulting estimate is bigger than $- 0.8\, \beta_1^2$.

(iii) Finally, we conclude that the estimate on the exponent obtained so far is approximately $- 0.8\, \beta_1^2$. Comparing this with our
estimate $- \beta_1^2$, it is clear that there is no much room left for an improvement by part (2) of Theorem 4.1.

\medskip

We now study the general criteria and estimates of $\lambda_1(L)$
and $\lambda_0^{\pm}(\theta)$ (see (4.17) below for definitions)
in dimension one. For this, we need more notation.

Fix an arbitrary reference point $\theta \in \Bbb R$, not
necessarily a root of $b(x)=-u'(x)$. Let $\Bbb R_{\theta
}^+=(\theta , \infty )$, $\Bbb R_{\theta }^-=(- \infty, \theta )$,
${\overline {\Bbb R}}_{\theta }^{\,+}=[\theta , \infty )$, and
${\overline {\Bbb R}}_{\theta }^{\,-}=(- \infty, \theta]$. Recall
that $C_{\theta }(x)=\int_{\theta }^x b/a$. Define
$$\align
&\text{\hskip-2em}\varphi_\theta^{} (x)=\int_{\theta }^x
e^{-C},\qquad \delta_{\theta }^{\pm}=\sup_{x\in \Bbb R_{\theta
}^{\pm}}
\varphi (x)\int_x^{\pm \infty } \frac{e^C}{a}, \\
&\text{\hskip-2em}{\scr F}_{I\theta }^{\pm}=\Big\{f\in
C\big({\overline {\Bbb R}}_{\theta }^{\pm}\big) \cap C^1({\Bbb
R_{\theta }^{\pm}}): f(\theta )=0,\;
f'\big|_{\Bbb R_{\theta }^{\pm}}>0\Big\},\\
&\text{\hskip-2em}{\scr F}_{I\!I\theta }^{\pm}=\Big\{f\in
C\big({\overline {\Bbb R}}_{\theta }^{\pm}\big):
f(\theta )=0, \; (\pm f)\big|_{\Bbb R_{\theta }^{\pm}}>0\Big\},\\
&\text{\hskip-2em} I_{\theta
}^{\pm}(f)(x)=\frac{e^{-C(x)}}{f'(x)}\int_{x}^{\pm \infty } \frac{f
e^C}{a},
\qquad \pm(x-\theta )\ge 0,\; \; f\in {\scr F}_{I\theta }^{\pm},\\
&\text{\hskip-2em} I\!I_{\theta }^{\pm}(f)(x)=\frac{1}{f(x)}\int_{\theta }^{\pm
\infty } \! \varphi_\theta^{} (x\wedge \cdot)\, \frac{f e^C}{a}
=\frac{1}{f(x)}\int_{\theta }^x e^{-C(y)}\d y\int_{y}^{\pm \infty
}\frac{f e^C}{a},\\
&\text{\hskip-2em} \qquad \qquad
\pm(x-\theta )\ge 0,\; \; f\in {\scr F}_{I\!I\theta }^{\pm},\\
&\text{\hskip-2em}\lambda_0^{\pm} (\theta)=\inf\Big\{D(f): f|_{\Bbb R\setminus \Bbb
R_{\theta }^{\pm}}=0,\, f\in C\big({\overline {\Bbb R}}_{\theta
}^{\pm}\big) \cap C^1({\Bbb R_{\theta }^{\pm}}),\, \mu(f^2)=1
\Big\}. \text{\hskip-2em}\tag 4.17
\endalign$$

\proclaim{Theorem 4.4} The comparison of $\lambda_1(L)$ and
$\lambda_0^{\pm} (\theta)$ and their estimates are given as
follows. \roster \item $\inf_{\theta \in \Bbb R} \big[
\lambda_0^+(\theta)\vee \lambda_0^-(\theta)\big] \ge \lambda_1(L)
\ge \sup_{\theta \in \Bbb R} \big[\lambda_0^+(\theta)\wedge
\lambda_0^-(\theta) \big]$. In particular, $\lambda_1(L) =
\lambda_0^+(\bar\theta)$, where $\bar \theta $ is the solution to
the equation $\lambda_0^+(\theta)=\lambda_0^-(\theta)$, $\theta
\in [-\infty , \infty ]$. \item If $m$ is the medium of $\mu$,
then $ 2\big[\lambda_0^+(m)\wedge \lambda_0^-(m)\big]\ge
\lambda_1(L)\ge  \lambda_0^+(m)\wedge \lambda_0^-(m).$ \item
$\lambda_0^{\pm} (\theta)\ge \sup_{f\in {\scr F}_{I\!I\theta
}^{\pm}} \inf_{x\in \Bbb R_{\theta }^{\pm}} I\!I_{\theta
}^{\pm}(f)(x)^{-1} \ge \sup_{f\in {\scr F}_{I\theta }^{\pm}}
\inf_{x\in \Bbb R_{\theta }^{\pm}} I_{\theta }^{\pm}(f)(x)^{-1}$.
Moreover, the sign of equalities hold whenever both $a$ and $b$
are continuous. \item $\big(\delta_{\theta }^{\pm}\big)^{-1} \ge
\lambda_0^{\pm} (\theta)\ge \big(4\delta_{\theta
}^{\pm}\big)^{-1}$.
\endroster
\endproclaim

\demo{Proof} The first assertion of part (1) is just [18; Theorem
3.3]. The lower bound in part (2) follows from the one of part (1).
As remarked above [18; Theorem 3.3], from the proof of [18; Theorem
3.1], it follows that
$$\lambda_1(L)\le \inf_{\theta\in \Bbb R}
\big[\lambda_0^+(\theta)\mu(\theta, \infty)\big]\wedge
\big[\lambda_0^-(\theta) \mu(-\infty, \theta)\big].$$ Hence, the
upper bound in part (2) follows immediately. The variational
formulas for the lower bounds given in part (3) is a copy of [19;
Theorem 1.1]. In which, the corresponding variational formulas for
the upper bounds are also presented, but omitted here. Part (4) was
proven in [18; Theorem 1.1]. From these quoted papers, one can find
some more sharper estimates and further references.

It remains to prove the second assertion of part (1). For this, it
suffices to show that $\lambda_0^{\pm}(\theta )$ is continuous in
$\theta $. By symmetry, it is enough to prove that
$\lambda_0^+(\theta )$ is continuous in $\theta $. Let $\theta_1<
\theta_2<\infty $. Clearly, $\lambda_0^+(\theta_1 )
<\lambda_0^+(\theta_2 )$. Given $\varepsilon \in (0, 1)$, choose
$f=f_{\varepsilon}\in C^1(\theta_1, \infty)\cap C[\theta_1, \infty
)$ such that $f(\theta_1)=0$, $ \int_{\theta_1}^\infty {f}^2 \d
\mu =1$ and $A-\varepsilon \le \lambda_0^+(\theta_1 )$, where
$A=A_{\varepsilon}=\int_{\theta_1}^\infty {f'}^2 \d \mu.$ By the
continuity of $f$, when $\theta_2-\theta_1>0$ is sufficient small,
we have
$$\align
\bigg|f(\theta_2)^2\int_{\theta_2}^\infty \d \mu &- 2 f(\theta_2)
\int_{\theta_2}^\infty f \d \mu- \int_{\theta_1}^{\theta_2}f^2 \d
\mu \bigg|\\
&\le f(\theta_2)^2\int_{\theta_1}^\infty \d\mu + 2
f(\theta_2)+\int_{\theta_1}^{\theta_2}f^2\d\mu\\
&<\varepsilon .\endalign
$$
Then $\int_{\theta_2}^\infty [f- f(\theta_2)]^2 \d \mu> 1-\varepsilon$
and furthermore
$$
\lambda_0^+(\theta_2 ) \le \int_{\theta_2}^\infty {f'}^2 \d \mu
\bigg/ \int_{\theta_2}^\infty [f- f(\theta_2)]^2 \d \mu
\le
\frac{A}{1-\varepsilon }\le \frac{\lambda_0^+(\theta_1
)+\varepsilon }{1-\varepsilon}.
$$  Since $\varepsilon $ can be arbitrarily
small, we obtain the required assertion. \qed\enddemo

As an illustration of the applications of Theorem 4.4, we discuss
Examples 4.2 and 4.3 again.

\proclaim{Example 4.5} Everything in premise is the same as in
Example 4.2. We have $$\frac{2 \alpha}{\delta}\ge
\lambda_1(L_{\alpha, \beta})\ge \frac {\alpha}{4 \delta},$$ where
$$\delta=\sup_{x>0}\int_0^x e^{y^2}\d y \int_x^\infty e^{-y^2} \d
y\approx 0.239405.$$
\endproclaim

\demo{Proof} First, we have the root $\theta =-\beta /(2\alpha )$
of $u'(x)$, it is also the medium of the measure. Next,
$$C_{\theta }(x)=-\alpha (x -\theta )^2,\quad
\varphi_{\theta }^{} (x)=\int_0^{x-\theta }e^{\alpha y^2}\d y,
\quad \int_{x}^\infty\! e^{-\alpha (y -\theta)^2}\d y =
\int_{x-\theta }^\infty e^{-\alpha y^2}\d y.
$$
Hence $\delta_{\theta }^+={\delta/\alpha}.$ By symmetry, we also
have $\delta_{\theta }^-={\delta/\alpha}.$ The assertion now
follows from parts (2) and (4) of Theorem 4.4. \qed\enddemo

\proclaim{Example 4.6} Everything in premise is the same as in
Example 4.3. We have \roster \item $\lim_{|\beta_{2}|\to
\infty}\lambda_1(L_{\beta_1,\, \beta_2})=\infty$.\item
 For $\beta_1\ge 0$, we have
$$
\lambda_1(L_{\beta_1,\, 0})\le
{4 e^{14}} \exp\bigg[-\frac 1 4 \bz_1^2+2\log(1+\bz_1)\bigg].
$$\endroster
\endproclaim

\demo{Proof} By symmetry of $u(x)$ in $x$, one may assume that
$\beta_2\ge 0$. Let $\theta$ be a real root of $u'(x)$. Clearly,
$\lim_{\beta_2\to\infty}\theta=-\infty$. Moreover, $u(x)-u(\theta
)=(x-\theta )^2 \big[(x-\theta )^2+4\theta (x-\theta
)+6\theta^2-\beta_1 \big]$. Hence
$$\align
\int_{\theta }^x& e^{u(y)}\d y \int_x^\infty e^{-u(z)}\d z \\
&=\int_{\theta }^x e^{u(y)-u(\theta )}\d y \int_x^\infty e^{-u(z)+u(\theta )}\d z \\
&=\int_0^{x-\theta }e^{y^2(y^2+ 4 \theta y+6 \theta^2-\beta_1)}\d
y
\int_{x-\theta }^\infty e^{- z^2(z^2+ 4 \theta z+6 \theta^2-\beta_1)}\d z \\
&=\int_0^{x-\theta }\!\! \d y \int_{x-\theta }^\infty \!
\exp\bigg[\!-\! (z^2\!-\!y^2)\bigg(z^2\!+\!y^2\!+\! 4 \theta
\bigg(z\!+\!\frac{y^2}{z+y} \bigg)
\!+\!6\theta^2\!-\!\beta_1\bigg)\bigg]\d z.\\
& \tag 4.18
\endalign$$

(a) We now prove the first assertion. It says that the parameter
$\beta_2$ plays a role for $\lambda_1(L_{\beta_1,\, \beta_2})$, in
contrast with Example 4.5. For $x\ge \theta$, by (4.18), we have
$$\align
&\int_{\theta }^x e^{u(y)}\d y \int_x^\infty e^{-u(z)}\d z \\
&=\int_{0}^{x-\theta}\d y \int_{x-\theta}^\infty\d z
\exp\bigg[-(z^2-y^2)\bigg[(z+2\theta)^2+\bigg(y+\frac{2\theta
y}{z+y}\bigg)^2\\
&\text{\hskip 12em}-4\theta^2\bigg(\frac{y}{z+y}\bigg)^2+2\theta^2-\beta_1\bigg]\bigg]\\
&\le \int_{0}^{x-\theta}\d y \int_{x-\theta}^\infty\d z
\exp\bigg[-(z^2-y^2)\bigg[-4\theta^2\bigg(\frac{y}{z+y}\bigg)^2+2\theta^2-\beta_1\bigg]\bigg].
\endalign$$
Since $z\ge y\ge 0$, we have $y/(z+y)\le 1/2$. The right-hand side
is controlled by
$$\int_{0}^{x-\theta}\d y\int_{x-\theta}^\infty e^{-(z^2-y^2)(\theta^2-\beta_1)}\d z,
\qquad x\ge \theta. \tag 4.19$$ We now use Conte's estimate (cf.
[20]):
$$x\bigg(1+\frac{x^2}{12}\bigg)e^{-3x^2/4}<
e^{-x^2}\int_0^x e^{y^2}\le \frac{\pi^2}{8x}(1-e^{-x^2}),\qquad
x>0$$ and Gautschi's estimate (cf. [21]):
$$\align
\frac 1 2 \Big[(x^p+2)^{1/p}-x\Big] < & e^{x^p}\int_x^\infty
e^{-y^p}dy\le C_p
\bigg[\bigg(x^p+\frac{1}{C_p}\bigg)^{1/p} - x\bigg],\qqd x\ge 0,\\
&C_p:=\ggz\big(1+{1}/{p}\big)^{p/(p-1)},\qd p>1;\qqd C_2=\pi/4.
\endalign$$
Thus,
$$\align
\int_{0}^{x}e^{c y^2}\d y \int_{x}^\infty e^{-c z^2}\d z &\le
\frac{\pi^2}{8c\sqrt{c}\,x}\big(1-e^{-c
x^2}\big)\cdot\frac{\pi}{4} \bigg(\sqrt{c x^2
+\frac{4}{\pi}}-\sqrt{c}\,x \bigg)\\
&\le \frac{\pi^2}{8c\sqrt{c}\,x}
\sqrt{\frac{\pi}{4}}\,\big(1-e^{-c x^2}\big),\qquad x\ge 0.
\endalign$$
Noting that $\big(1-e^{-c x^2}\big)/x\le c x \le c$ for all $x\in
(0, 1]$ and $\big(1-e^{-c x^2}\big)/x\le {1}/{x}\le 1$ for all $x\ge
1$, we obtain
$$
\int_{0}^{x} e^{c y^2}\d y \int_{x}^\infty e^{-c z^2}\d z
\le
\frac{\pi^{5/2}}{16\sqrt{c}}, \qquad x\ge 0,\; c\ge 1.$$ Therefore
$$\align
\delta_{\theta}^+&= \sup_{x>\theta}\int_{\theta }^x e^{u(y)}\d y
\int_x^\infty e^{-u(z)}\d z \\
& \le \sup_{x>0} \int_{0}^{x}\d
y\int_{x}^\infty
e^{-(z^2-y^2)(\theta^2-\beta_1)}\d z \\
&\le \frac{\pi^{5/2}}{16\sqrt{\theta^2-\beta_1}}\to 0
\qquad\text{as }\; \theta\to -\infty.\endalign
$$

For $\delta_\theta^-$, the proof is similar. As an analogue of
(4.18), we have
$$\align
&\int_x^{\theta } e^{u(y)}\d y \int_{-\infty}^x e^{-u(z)}\d z \\
&=\int_{x-\theta}^0\d y \int_{-\infty}^{x-\theta}\d z
\exp\bigg[-(z^2-y^2)\bigg[(z+2\theta)^2+\bigg(y+\frac{2\theta
y}{z+y}\bigg)^2\\
&\text{\hskip12em}-4\theta^2\bigg(\frac{y}{z+y}\bigg)^2+2\theta^2-\beta_1\bigg]\bigg].
\endalign$$
Since $z\le y\le 0$, we have $|y/(z+y)|\le 1/2$, we obtain
$$\int_x^{\theta } e^{u(y)}\d y \int_{-\infty}^x e^{-u(z)}\d z
\le \int_{x-\theta}^0 \d y \int_{-\infty}^{x-\theta}
e^{-(z^2-y^2)(\theta^2-\beta_1)}\d z,\qquad x\le \theta.$$ We have
thus returned to (4.19).

Now, the first assertion follows from parts (1) and (4) of Theorem
4.4.

(b) For the upper bound in part (2), since $\beta_2=0$, we
have $\theta=0$. We need to show that
$$\sup_{x>0} \int_0^x  e^{y^4-\beta_1 y^2} \d y \int_x^\infty e^{-z^4 +
\beta_1 z^2} \d z\ge \frac{1}{4 e^{14}} \exp\bigg[\frac 1 4 \bz_1^2-2\log(1+\bz_1)\bigg].$$
Since
$$\align
\int_0^x & e^{y^4-\beta_1 y^2} \d y \int_x^\infty e^{-z^4 +
\beta_1 z^2} \d z\\
&= \frac 1 4 \int_{-\bz_1/2}^{x^2-\bz_1/2} \frac{ e^{y^2}}{\sqrt{y+\bz_1/2}} \d y
\int_{x^2-\bz_1/2}^\infty \frac{e^{-z^2}}{\sqrt{z+\bz_1/2}} \d z\\
&> \frac 1 4 \int_{-\bz_1/2}^{x^2-\bz_1/2} \frac{ e^{y^2}}{\sqrt{y+\bz_1/2}} \d y
\int_{x^2-\bz_1/2}^{\bz_1 /2} \frac{e^{-z^2}}{\sqrt{z+\bz_1/2}} \d z,\endalign$$
when $\bz_1\ge 1$, we have
$$\align
 \int_{-\bz_1/2}^{1-\bz_1/2}& \frac{ e^{y^2}}{\sqrt{y+\bz_1/2}} \d y
\int_{1-\bz_1/2}^{\bz_1 /2} \frac{e^{-z^2}}{\sqrt{z+\bz_1/2}} \d z\\
&\ge \frac{1}{\bz_1} \int_{-\bz_1/2}^{1-\bz_1/2}  e^{y^2} \d y
\int_{1-\bz_1/2}^{\bz_1 /2} e^{-z^2} \d z.\endalign$$
It suffices to show that
$$\frac{1}{\bz_1} \int_{-\bz_1/2}^{1-\bz_1/2}  e^{y^2} \d y
\int_{1-\bz_1/2}^{\bz_1 /2} e^{-z^2} \d z\ge \frac{1}{e^{14}} \exp\bigg[\frac 1 4 \bz_1^2-2\log(1+\bz_1)\bigg],$$
or
$$\int_{-\bz_1/2}^{1-\bz_1/2}  e^{y^2} \d y
\int_{1-\bz_1/2}^{\bz_1 /2} e^{-z^2} \d z\ge \exp\bigg[\frac 1 4 \bz_1^2-\log(1+\bz_1)-14\bigg].$$
Since
$$\align
\int_{1-\bz_1/2}^{\bz_1 /2} e^{-z^2} \d z&\to \int_{-\infty}^{\infty} e^{-z^2} \d z<\infty,\\
\int_{-\bz_1/2}^{1-\bz_1/2}  e^{y^2} \d y&=\int_{\bz_1/2-1}^{\bz_1/2}  e^{y^2} \d y\ge \exp\bigg[\bigg(\frac{\bz_1}{2}-1\bigg)^2\bigg]\to\infty,\\
\frac{\int_{-\bz_1/2}^{1-\bz_1/2}  e^{y^2} \d y}{\exp[\bz_1^2/4-\log \bz_1]}
&\sim\frac{\exp[\bz_1^2/4]-\exp[(1-\bz_1/2)^2]}{\exp[\bz_1^2/4]}\\
&\sim 1-e^{1-\bz_1}\\
&\sim 1\qd\text{as }\bz_1\to\infty,
\endalign$$
it is easy to check first that
$$\log\bigg[ \int_{-\bz_1/2}^{1-\bz_1/2}  e^{y^2} \d y
\int_{1-\bz_1/2}^{\bz_1 /2} e^{-z^2} \d z\bigg]\ge \frac 1 4 \bz_1^2-\log(1+\bz_1)-14$$
for $\bz_1\ge 1$ and then the required assertion for $\bz_1\ge 0$ by using mathematical softwares.
 \qed \enddemo

Before moving further, let us study the lower bounds of
$\inf_{\beta_2\ge 0}\lambda_1(L_{\beta_1,\, \beta_2})$ in terms of $\delta_{\theta}^{\pm}$.
For this, we return to (4.18). Because
$$
4 \theta \bigg(z+\frac{y^2}{z+y} \bigg) +6\theta^2 =6\bigg[\theta+
\frac{1}{3}\bigg(z+\frac{y^2}{z+y}\bigg)
\bigg]^2-\frac{2}{3}\bigg(z+\frac{y^2}{z+y} \bigg)^2\ge
-\frac{2}{3}\big(z+ y/2 \big)^2,
$$
and so
$$\align
 z^2+y^2+ 4 \theta \bigg(z+\frac{y^2}{z+y} \bigg)
+6\theta^2-\beta_1 &\ge z^2+y^2-\frac{2}{3}\big(z+ y/2 \big)^2-
\beta_1\\
&=\frac{1}{6} (2
z^2 - 4 zy + 5y^2 ) -\beta_1\\
&\ge \frac {1}{6} (z^2+y^2)-\beta_1,
\endalign
$$
it follows that
$$\align
\int_0^{x} \d y& \int_{x}^\infty \exp\bigg[-
(z^2-y^2)\bigg(z^2+y^2+ 4 \theta \bigg(z+\frac{y^2}{z+y} \bigg)
+6\theta^2-\beta_1\bigg)\bigg]\d z\\
 &\le \int_0^x \d y
\int_x^\infty \d z\, e^{-(z^2-y^2)((z^2+y^2)/6
-\beta_1)}\\
&=\int_0^x \d y\,  e^{y^4/6-\beta_1 y^2} \int_x^\infty e^{-z^4/6+
\beta_1 z^2} \d z. \tag 4.20
\endalign$$
Combining (4.18) with (4.20), we obtain
$$\delta_{\theta}^+\le \sup_{x>0}\int_0^x \d y \, e^{y^4/6-\beta_1 y^2} \int_x^\infty e^{-z^4/6+
\beta_1 z^2} \d z.
 $$
The same upper bound holds for $\delta_{\theta}^-$. By parts (1)
and (4) of Theorem 4.4, we obtain a lower estimate of
$\inf_{\beta_2\ge 0}\lambda_1(L_{\beta_1,\, \beta_2})$. However,
the resulting bound is smaller than those given in Example 4.3.

We mention that the lower bound given in Example 4.3 may still be
improved by applying part (3) of Theorem 4.3 to the test functions
$f_{\pm}$ constructed in the proof of Theorem 4.1. This observation
is due to [22]. The proof is quite easy. Let for instance
$$-\sup_{x\in (\theta, \infty)}\frac{a f_+'' + b f_+'}{f_+}(x)\ge \delta>0.$$
Then $f_+\le -(af_+''+bf_+)/\delta$. Noting that $\big(e^C
f_+'\big)'= e^C (af_+''+bf_+)/a$, we obtain
$$\align
I_{\theta}^+(f_+)(x)&=\frac{e^{-C(x)}}{f_+'(x)}\int_x^\infty
\frac{f_+ e^C}{a}\\
& \le
\frac{1}{\delta}\frac{e^{-C(x)}}{f_+'(x)}\int_x^\infty
\bigg(-\frac{af_+''+bf_+}{a}\bigg)e^C\\
& = \frac{1}{\delta}\frac{e^{-C(x)}}{f_+'(x)}\int_x^\infty
\big(-e^C f_+'\big)'\\
&\le \frac{1}{\delta}\frac{e^{-C(x)}}{f_+'(x)}
e^{C(x)} f_+'(x)\\
&=\frac{1}{\delta},\qquad x>\theta. \endalign$$
 Alternatively, one may
apply the approximation procedure given in [19] to improve the lower
bound. However, all the computations are quite complicated, and so
we do not want to go further along this line.

We remark that the process in Example 4.6 (Example 4.3) possesses much stronger
ergodic properties.

\proclaim{Proposition 4.7}
The processes corresponding to Example 4.3
is not only exponentially ergodic but also
strongly ergodic. It has the empty essential spectrum. It satisfies
the logarithmic Sobolev inequality but not the Nash (Sobolev)
inequality. \endproclaim

\demo{\prf} One may use the criteria given in [13; \S 5.4] to
justify these assertions. For the reader's convenience, here we
mention three criteria as follows. By the symmetry, we need only to
write down the conditions on the half-line $[0, \infty)$.

{\it Logarithmic Sobolev inequality}:
$$\sup_{x>0}\bigg(\int_x^\infty e^{-u}\bigg)\bigg(\log\int_x^\infty e^{-u} \bigg)
\int_0^x e^u<\infty.$$

{\it Strong ergodicity}:
$$\int_0^\infty\d x e^{u(x)}\int_x^\infty e^{-u}<\infty.$$

{\it Nash $($Sobolev$)$ inequality}:
$$\sup_{x>0}\bigg(\int_x^\infty e^{-u}\bigg)^{1-2/\nz}
\int_0^x e^u<\infty,\qqd \nz>2.$$

The second condition holds since
$$\frac{\int_x^\infty e^{-u}}{x^{-2}e^{-u}}
\sim \frac{1}{2 x^{-3}+x^{-2}u'}\sim \frac{x^3}{x u'}\to 0, \qqd x\to\infty.$$
However, replacing $x^{-2}$ with $x^{-1}$ at the beginning, the same proof shows that the
standard Ornstein-Uhlenbeck process is not strongly ergodic.
 For the third condition, note that
$\int_x^\infty e^{-u}$ and $\int_0^x e^u$ have the leading order
$e^{-u}$ and $e^u$ respectively. Hence the leading order of
$$\bigg(\int_x^\infty e^{-u}\bigg)^{1-2/\nz}
\int_0^x e^u$$ is $e^{2 u/\nz}\to\infty$ as $x\to\infty.$ Similarly,
one can check the first condition. Alternatively, to see that the
logarithmic Sobolev inequality holds, simply use the fact that
$\lim_{|x|\to \infty}u''(x)>0$ (see [23]). We will come back to this
point in Example 5.3. Finally, the logarithmic Sobolev inequality
implies the essential spectrum to be empty. \qed\enddemo

Finally, we
study a perturbation of $\lambda_1(L)$.

\proclaim{Proposition 4.8} Let $a(x)\equiv 1$ and assume that $\delta_{\theta
}^{\pm}<\infty $ for some $\theta \in \Bbb R$. Next, let $h$ satisfy
$\int_{\Bbb R} e^{C+h}<\infty $. Define $\delta_{\theta }^{\pm}(h)
=\sup_{x\in {\Bbb R}_{\theta }^{\pm}} \int_{\theta }^x
e^{-C-h}\int_x^{\pm\infty } e^{C+h}$. If there exist constants
$K_1^{\pm}, \ldots, K_4^{\pm}$ such that
$$\align
&\pm \int_x^{\pm\infty} e^C\le K_1^{\pm} e^{C(x)},\qquad \pm(x-\theta )\ge 0, \tag 4.21\\
&\pm \int_{\theta}^{x} e^{-C}\le K_2^{\pm} e^{-C(x)},\qquad \pm(x-\theta )\ge 0, \tag 4.22\\
&\pm \int_x^{\pm\infty} e^C\big|e^h-1\big|\le K_3^{\pm} e^{C(x)},\qquad \pm(x-\theta )\ge 0, \tag 4.23\\
&\pm \int_{\theta}^{x} e^{-C}\big|e^{-h}-1\big|\le K_4^{\pm}
e^{-C(x)},\qquad \pm(x-\theta )\ge 0, \tag 4.24
\endalign$$
then
$$ \delta_{\theta }^{\pm}(h)\le
\delta_{\theta }^{\pm}+  K_2^{\pm} K_3^{\pm} + K_1^{\pm} K_4^{\pm}
+ K_4^{\pm} K_3^{\pm}< \infty. $$
\endproclaim

\demo{Proof} Here, we consider $\delta_{\theta }^+(h)$ only. As in
[5], we have
$$\align
&\int_{\theta }^x e^{-C-h}\int_x^\infty e^{C+h} \\
&= \bigg[\int_{\theta }^x e^{-C}+\int_{\theta }^x
e^{-C}\big(e^{-h}-1\big)\bigg]
\cdot \bigg[\int_x^\infty  e^{C}+\int_x^\infty e^C \big(e^{h}-1\big)\bigg]\\
&= \int_{\theta }^x e^{-C} \int_x^\infty  e^{C}
+\int_{\theta}^x e^{-C} \int_x^\infty  e^{C} \big(e^{h}-1\big)\\
&\quad + \int_{\theta}^x e^{-C}\big(e^{-h}-1\big)  \int_x^\infty
e^{C} +\int_{\theta}^x e^{-C}\big(e^{-h}-1\big) \int_x^\infty
e^{C}
\big(e^{h}-1\big) \\
&\le \delta_{\theta }^{+}+ K_2^+ K_3^+ + K_1^+ K_4^+ + K_4^+ K_3^+
< \infty. \qed
\endalign$$
The above result is a revised version of [5; Theorem 3.4], where
instead of (4.23) and (4.24), the conditions \roster
\item"(i)" ${C(x)}$
is strictly uniformly convex up to a bounded function
\item"(ii)"
$\int_{\Bbb R} \big(e^{|h|}-1\big)<\infty$\endroster
 are employed.
It is easy to check that these conditions together are stronger
than (4.23) and (4.24). Clearly, under (4.21) and (4.22),
conditions (4.23) and (4.24) are automatic for bounded $h$, for
which, the condition (ii) here 
may fail.
\enddemo

\proclaim{Example 4.9} Let $a(x)\equiv 1$ and $C_{\bz}(x)=-x^4+\bz x^2$. Then
$\lz_1(L_{\bz})>0$ for all $\bz\in {\Bbb R}$.
\endproclaim

\demo{\prf} The case of $\bz<0$ is easy since $-C_{\bz}$ is convex.
Hence we assume that $\bz\ge 0$. Then $-C_{\bz}$ is convex for large enough $x$
and so the conclusion is known. Here we check it by using
Proposition 4.8. Take $C(x)=-x^4$ and regard $h(x)=\bz x^2$ as a
perturbation of $C(x)$. Clearly, $\int_{\Bbb R}
\big(e^{|h|}-1\big)=\infty$. Set $\uz=0$.

First, by Gautschi's estimate, we have
$$e^{-C(x)}\int_x^\infty e^{C}=e^{x^4}\int_x^\infty e^{-y^4}\d y \le C_4 \bigg[\bigg(x^4+\frac{1}{C_4}\bigg)^{1/4}-x\bigg]
\le \ggz \bigg(\frac 5 4\bigg)\approx 0.9064$$ for all $x>0$. Next, we
have
$$\align
e^{C(x)}\int_x^\infty
e^{-C}\big|e^{-h}-1\big|&=e^{-x^4}\int_x^\infty e^{y^4}\big|e^{-\bz
y^2}-1\big|\d y\\
& \le e^{C(x)}\int_x^\infty e^{-C}\\
&= e^{-x^4}\int_x^\infty e^{y^4}\d
y \\
&<0.6,\qqd x>0.\endalign$$ Moreover,
$$e^{-C(x)}\int_x^\infty e^{C}\big|e^{-h}-1\big|< e^{-C(x)}\int_x^\infty e^{C+h}\le e^{\bz (0.7\,\bz + 26)},\qqd x>0.$$
By symmetry, the same estimates hold on $(-\infty, 0]$. Now, by Proposition 4.8 and Theorem 4.4, it
follows that the leading order
of the lower estimate of $\lz_1(L_{\bz})$ is $\exp[-0.7\, \bz^2]$ which is not far away from the optimal
one: $\exp[-\bz^2/4]$.
\qed\enddemo

\subhead{5. Logarithmic Sobolev inequality}
\endsubhead

We begin this section with a result taken from [23; Corollary 1.4].

\proclaim{Lemma 5.1} Let $L=\Delta -\langle \nabla U, \nabla\rangle$
in $\Bbb R^n$ and define $\gz(r)\!=\!\inf\limits_{|x|\ge
r}\lambda_{\min}(\hes (U)(x))$. If $\sup_{r\ge 0}\gz(r)$ $>0$, then
we have
$$\sz(L)\ge \frac{2 e}{a_0^2}\exp\bigg[- \int_0^{a_0^{}} r \gz (r) \d r\bigg]>0, $$
where $a_0^{}>0$ is the unique solution to the equation $\int_0^a \gz (r) \d r=2/a$.
\endproclaim

This lemma says that the logarithmic Sobolev constant is positive
whenever so is \allowlinebreak$\lambda_{\min}(\hes (U)(x))$ at infinity.
Unfortunately, as shown by Example 2.5, our models do not satisfy
this condition even in the two-dimensional case. Hence, we justify
the power of the estimate provided by the lemma only in dimensional
one (compare with the criterion for the inequality, see for instance
[13; Theorem 7.4]).

\proclaim{Example 5.2} For
Example 4.2, we have $\lz_1(L_{\az,\,\beta})\ge\sz(L_{\az,\,\beta})\ge 2\az$ which are exact.
\endproclaim

\demo{Proof} Because $u(x)=\az x^2+\bz x$, we have $ u''(x)=2\az$ and so
$$\gz(r)=\inf\limits_{|x|\ge r}u''(x)= 2\az.$$ Next, since $\int_0^a \gz(r)\d
r= 2\az\, a.$ The unique solution to the equation
$$\int_0^a \gz(r)\d r=\frac{2}{a}$$ is $a_0^2=1/\az.$ Noticing that
$\int_0^a r \gz(r)\d r= {\az}\,a^2,$
by Lemma 5.1, we obtain
$$\sz(L_{\az,\,\beta})\ge \frac{2 e}{a_0^2}
\exp\big[-\az\, a_0^2\big]=2\az.$$ This is clearly exact since the
well-known fact $\lz_1(L_{\az,\,\beta})\ge\sz(L_{\az,\,\beta})$ (cf.
[13; Theorem 8.7]) and Example 4.2.\qed\enddemo

\proclaim{Example 5.3} For
Example 4.3, we have
$$\align
\inf_{\bz_2}\lz_1(L_{\bz_1,\,\beta_2})&\ge\inf_{\bz_2}\sz(L_{\bz_1,\,\beta_2})\ge \frac{\sqrt{\beta_1^2+8}- \beta_1}{\sqrt{e}}
\exp\bigg[-\frac{1}{8}\beta_1\Big(
\beta_1+\sqrt{\beta_1^2+8}\Big)\bigg]\\
&\ge {\cases -2\beta_1+\dfrac{2}{\sqrt{e/2}-\bz_1},
\quad &\text{if } \beta_1<0\\
2\sqrt{2/e},\quad &\text{if } \beta_1=0\\
\dfrac{1}{\sqrt{e/8}+\bz_1}\,\exp\big[-\beta_1^2/4\big],\quad &\text{if }
\beta_1>0
\endcases}\endalign$$
\endproclaim

\demo{Proof} Because $u(x)=x^4-\bz_1 x^2+ \bz_2 x$, we have $
u''(x)=12 x^2-2 \bz_1$ and $\gz(r)=\inf_{|x|\ge r}u''(x)= 12 r^2-2
\bz_1$. Next, since $\int_0^a \gz(r)\d r=4 a^3-2\bz_1 a,$ the
solution to the equation $\int_0^a \gz(r)\d r=2/a$ is as follows
$$a_0^2=\frac{\bz_1+\sqrt{\bz_1^2+8}}{4}.$$
Next, since
$$\int_0^a r \gz(r)\d r=a^2(3 a^2-\beta_1),$$
by Lemma 5.1, we obtain
$$\align
\sz(L_{\beta_1,\,\beta_2})&\ge \frac{2 e}{a_0^2} \exp\big[-a_0^2(3
a_0^2-\beta_1)\big]\\
&= \frac{\sqrt{\beta_1^2+8}- \beta_1}{\sqrt{e}}
\exp\bigg[-\frac{1}{8}\beta_1\Big(
\beta_1+\sqrt{\beta_1^2+8}\Big)\bigg].\qed\endalign$$ \enddemo

Note that in the case of $\bz_1<0$, the Bakry-Emery criterion (cf.
(2.9)) is available and gives us the lower bound $-2\bz_1$ which is
smaller than the estimate above. Example 5.3 is somehow unexpect
since it improves Example 4.3 (In the special case that $\bz_2=0$,
they are coincided). The reason is due to the fact that only the
uniform estimate is treated in Example 4.3 and the linear term of
$U$ is ruled out in Lemma 5.1 (but the universal estimates depend on
the linear term, cf. [13; Theorem 7.4]). Otherwise, the two methods
may not be comparable in view of part (1) of Example 4.6. As
mentioned in [23; Example 1.12] that the bounded perturbations
should be carefully treated before applying \lmm\;5.1.

\demo{Proof of \prp\;$1.4$} Let $\bz_1\ge 0$. Note that
$$\sqrt{1+\frac{8}{\bz_1^2}}\le 1+\frac{4}{\bz_1^2}.$$
We have $\sqrt{\bz_1^2+8}\le \bz_1+4/\bz_1$. Hence
$$\exp\bigg[-\frac{1}{8}\beta_1\Big(
\beta_1+\sqrt{\beta_1^2+8}\Big)\bigg]\ge \frac{1}{\sqrt{e}}\exp\bigg[-\frac{1}{4}\beta_1^2\bigg].$$
Similarly, we have
$$\sqrt{\beta_1^2+8}- \beta_1=\frac{8}{\sqrt{\beta_1^2+8}+\beta_1}\ge \frac{4}{\bz_1+2 \bz_1^{-1}}.$$
By Example 5.3, we obtain
$\inf_{\bz_2}\sz(L_{\bz_1,\,\beta_2})\ge \exp[-\bz_1^2/4-\log(1+\bz_1)]$ for $\bz_1\ge 2$. Combining this with
Example 4.6, we get the required assertion.
\qed\enddemo

\subhead{6. Continuous spin systems}
\endsubhead

We begin this section with the ergodicity of our models in the finite dimensions.
Consider the particle system on $\llz$ with periodic boundary. Then
the generator is
$$
L_{\llz}=\ddz + \langle b, \, \nabla\rangle
$$ where
$$b_i(x)= -u'(x_i)- 2 J \sum_{j\in N(i)}(x_i-x_j)$$
for some $u\in C^{\infty}(\Bbb R)$, constant $J$, and $N(i)$ is the nearest neighbors of $i$.
For simplicity, assume that $J\ge 0$, but it is not essential in this section. Recall that for the
coupling by reflection, the coupling operator $\overline L$ has the
coefficients
$${a(x, y)=\bigg(\matrix I & I -2 \bar u {\bar u}^*\\
I -2 \bar u {\bar u}^* & I\endmatrix}\bigg),\qquad
{b(x, y)=\bigg(\matrix b(x)\\
b(y)\endmatrix\bigg)},$$ where $\bar u=\bar u(x, y)=(x-y)/|x-y|$.
Furthermore, for $f\in C[0, \infty)\cap C^2(0, \infty)$, we have
$${\overline L}f(|x- y|)=
4 f''(|x-y|)+ \frac{\langle x-y,\, b(x)-b(y)\rangle}{|x-y|}\,
{f'(|x-y|)},\qqd x\ne y$$ (cf. [13; Theorem 2.30]). To illustrate
the idea, we restrict ourselves to the second model.

\proclaim{\thm\;6.1} Let $u(x_i)=x_i^4-\bz x_i^2$ for all $i\in\llz$. Then the process is
exponentially ergodic
for any finite $\llz$. Moreover, the coupling by reflection $(X_t, Y_t)$ gives
us
$${\overline{\Bbb E}}^{x, y}f(|X_t-Y_t|)\le f(|x-y|)e^{-\vz t},\qqd t\ge 0,$$
where
$$\align
&f(r)=\int_0^r e^{-C(s)}\d s \int_s^\infty e^C \sqrt{\fz},\qqd r>0,\\
&C(r)=-\frac{1}{16 |\llz|}r^4+\frac{\bz}{4} r^2,\qqd \fz(r)=\int_0^r e^{-C},\\
&\vz=\vz(\llz, \bz)=4\inf_{r>0} \frac{\sqrt{\fz (r)}}{f(r)}>0.
\endalign$$
\endproclaim

\demo{\prf} Because $u'(x_i)= 4 x_i^3 -2\bz x_i$ and
$$
b_i(x) =-4  x_i^3 +2\bz x_i-2 J \sum_{j\in N(i)}(x_i-x_j)
$$
for all $i$. Thus,
$$b_i(x)-b_i(y)=-4 \big(x_i^3-y_i^3\big) +2\bz (x_i-y_i)-2 J \sum_{j\in N(i)}(x_i-y_i-x_j+y_j).$$
Hence
$$\align
\langle x-y,\, b(x)-b(y)\rangle&=
-4 \sum_i(x_i-y_i)^2\big(x_i^2+x_iy_i+y_i^2\big)
+2\bz \sum_i (x_i-y_i)^2\\
&\qqd - J \sum_i\sum_{j\in N(i)}(x_i-y_i-x_j+y_j)^2\\
&\le - \sum_i(x_i-y_i)^4 +2\bz \sum_i (x_i-y_i)^2\\
&\le - |\llz|^{-1}|x-y|^4 +2\bz |x-y|^2,
\endalign
$$
where $|\llz|$ is the cardinality of $\llz$. It follows that
$$\frac{\langle x-y,\, b(x)-b(y)\rangle}{|x-y|}\le - \frac{1}{|\llz|}|x-y|^3 +2\bz |x-y|.$$
If we take $f(r)=r$, then for all $x\ne y$, we have
$${\overline L}f(|x- y|)=
\frac{\langle x-y,\, b(x)-b(y)\rangle}{|x-y|}\, {f'(|x-y|)}\le -
\bigg(\frac{1}{|\llz|}|x-y|^2 -2\bz \bigg)|x-y|.$$
This is not enough for the exponential convergence
 except in the case that $\bz<0$ for which we have
$\inf_{r>0}\big(r^2/|\llz| -2\bz \big)=-2\bz>0$. Due to this reason, we need a
much carefully designed $f$. Define the function $f$ as in the theorem,
then we have
$$f'(r)=e^{-C(r)} \int_r^\infty e^C \sqrt{\fz},\qqd f''=-\frac{1}{4} \gz f'-\sqrt{\fz}.$$
We obtain
$$4 f''+\gz f'=-4 \sqrt{\fz}\le -\vz f$$
with
$$\gz(r)=-\frac{1}{ |\llz|} r^3 +2\bz r,\qqd \vz=4\inf_{r>0} \frac{\sqrt{\fz (r)}}{f(r)}.$$
By the Cauchy mean value theorem, it follows that
$$\align
\inf_{r>0} \frac{\sqrt{\fz}}{f}&\ge \inf_{r>0} \frac{(\sqrt{\fz})'}{f'}= \frac 1 2\inf_{r>0}
\fz^{-1/2}\bigg/\int_r^\infty e^C \sqrt{\fz}\\
&\ge\frac 1 2\inf_{r>0}
\frac{(\fz^{-1/2})'}{- e^C \sqrt{\fz}}
=\frac 1 4 \bigg(\inf_{r>0} \frac{e^{-C}}{\fz}\bigg)^2>0.\endalign $$
Therefore we obtain $\vz>0$. This proves our second assertion.

The exponential ergodicity is easy to check by using the so called
``drift condition'' with test function $x\to |x|^2$, but this is not
enough to get a convergence rate. We now prove the exponential
ergodicity with respect to $f\circ |\cdot|$. Note that here we do
not assume that $f\circ |\cdot|$ is a distance. Otherwise, the
assertion follows from [17; Theorem 5.23]. We have proved in the
last paragraph that ${\overline{\Bbb E}}^{x, y}f(|X_t-Y_t|)$ is
continuous in $y$. Moreover
$${\overline{\Bbb E}}^{x, \mu_{U}^{}}f(|X_t-Y_t|)
=\int_{\Bbb R^{|\llz|}} \mu_{U}^{}(\d y) {\overline{\Bbb E}}^{x, y}f(|X_t-Y_t|)
\le e^{-\vz t} \int_{\Bbb R^{|\llz|}} \mu_{U}^{}(\d y) f(|x-y|),$$
where
$\mu_{U}^{}$ is the probability measure having density $e^{-U}/Z_U,$
corresponding to the potential
$$U(x)= \sum_{i\in \Lambda } u(x_i)
+ J\sum_{i\in \Lambda}\sum_{j\in N(i)}  (x_i- x_j)^2.$$
Because the left-hand side controls the Wasserstein distance, with respect
to the cost function $f\circ |\cdot|$, of the laws of the processes starting
from $x$ and $\mu_{U}^{}$ respectively, we obtain an exponential ergodicity provided
$$\int_{\Bbb R^{|\llz|}} \mu_{U}^{}(\d y) f(|x-y|)<\infty.$$
To check this, noting that
$$-U(x)\le \sum_{i\in \Lambda }\big(-x_i^4+\bz x_i^2\big)\le -\frac{1}{|\llz|}|x|^4+\bz |x|^2$$
and $f(|x-y|)\le f(|x|+|y|)$, it suffices to consider the radius part. That is,
$$\int_{0}^\infty f(r+z)\exp[-z^4/|\llz|+\bz z^2]\d z<\infty \qqd \text{for every } r\ge 0.$$
This can be done by using a comparison:
$$\align
\frac{f(r+z)\exp[-z^4/|\llz|+\bz z^2]}{z^{-2}}
& =\frac{f(r+z)}{z^{-2}\exp[z^4/|\llz|-\bz z^2]}\\
& \sim \frac{e^{-C(r+z)}\int_{r+z}^\infty e^C \sqrt{\varphi}}{[-2z^{-3}+z^{-2}(4 z^3/|\llz|-2\bz z)]\exp[z^4/|\llz|-\bz z^2]}\\
&\sim \frac{\int_{r+z}^\infty e^C \sqrt{\varphi}}{z\exp[z^4/|\llz|-\bz z^2+C(r+z)]}\\
&\sim 0\qqd \text{as } z\to\infty.
\endalign$$

Finally, by [17; Theorem 9.18] and its remark, we also have
$\lz_1(U, \llz, \bz)>0$. \qed\enddemo

Theorem 6.1 is meaningful since it works for all finite dimensions. Note that
$\vz(\llz, \bz)\to 0$ as $|\llz|\to \infty$, which is natural since the model
exhibits a phase transition. However, this result does not describe an ergodic region
in the infinite dimensional situation.

For the remainder of this section, we apply the results obtained in the previous
sections to some specific continuous spin systems. Denote by
$\langle ij\rangle $ the nearest bonds in $\Bbb Z^d,\; d\ge 1$. Set
$N(i)=\{j: j \text{ is the endpoint of an bond $\langle ij\rangle $}
\}$. Then, $|N(i)|:=\text{the cardinality of the set }N(i)$ $= 2 d$.
Consider the Hamiltonian $H(x)=J\sum_{\langle ij\rangle }  (x_i-
x_j)^2$, where $J\ge 0$ is a constant. For a finite set $\Lambda
\subset \Bbb Z^d$ (denoted by $\Lambda\Subset {\Bbb Z}^d$) and a
point $\omega \in {\Bbb R}^{\Bbb Z^d}$, define the
finite-dimensional conditional Gibbs distribution $\mu_U^{\Lambda,
\omega }$ as follows.
$$\mu_U^{\Lambda,\, \omega }(\d x_{\Lambda }^{})={e^{-U_{\Lambda }^\omega
(x_{\Lambda }^{})}}\, \d x_{\Lambda }^{}\big/ {Z_{\Lambda }^\omega
}, \tag 6.1$$ where $x_{\Lambda}=(x_i, i\in \Lambda)$, $
Z_{\Lambda }^\omega $ is the normalizing constant and
$$ U_{\Lambda }^\omega(x_{\Lambda }^{})=\sum_{i\in \Lambda } u(x_i)
+ J\sum_{\langle ij\rangle :\, i, j\in \Lambda }  (x_i- x_j)^2 +
J\sum_{i\in \Lambda ,\, j\in N(i)\setminus \Lambda } (x_i-
\omega_j)^2 \tag 6.2$$ for some function $u\in C^\infty (\Bbb R
)$, to be specified latterly. One can rewrite $U_{\Lambda
}^\omega$ as
$$U_{\Lambda }^\omega (x_{\Lambda }^{})= \sum_{i\in \Lambda } u(x_i)
+ J\sum_{i\in \Lambda}\sum_{j\in N(i)}  (x_i- z_j)^2, \tag 6.3 $$
where
$$z_j=
\cases x_j, &\quad \text{if } j\in \Lambda\\
\omega_j, &\quad \text{if } j\notin \Lambda.
\endcases$$
Correspondingly, we have an operator $L_{\Lambda }^\omega$ and a
Dirichlet form $D_{\Lambda }^\omega $ as follows.
$$ L_{\Lambda }^\omega =\Delta_{\Lambda }-\langle \nabla_{\Lambda}
U_{\Lambda }^\omega ,\, \nabla_{\Lambda }\rangle, \qquad
D_{\Lambda }^\omega (f)=\int_{\Bbb R^{|\Lambda |}}
|\nabla_{\Lambda } f|^2 \d \mu_U^{\Lambda,\,\omega}. \tag 6.4$$
Our purpose in this section is to estimate
$\lambda_1\big(L_{\Lambda }^\omega \big)= \lambda_1\big(U_{\Lambda
}^\omega \big)$. By (1.6), we have the simplest lower bound of the
marginal eigenvalues as follows.
$$\lambda_1^{x_{\Lambda\backslash i},\,\omega} \ge \inf_{x\in \Bbb R} u''(x) + 4 d J, \tag 6.5$$
where $x_{\Lambda\backslash i}=(x_j, j\in \Lambda \setminus \{i\})$. The function $C(x)$ defined in Section 4 becomes
$$\align
 C_{\Lambda }^{x_{\Lambda\backslash i},\,\omega }(x_i)&= -u(x_i)- J \sum_{j\in N(i)}
 (x_i-z_j)^2\\
&=-u(x_i)- 2dJx_i^2+2J \bigg(\sum_{j\in N(i)}z_j\bigg)x_i -J
\sum_{j\in N(i)}z_j^2,\\
& \qquad i\in \Lambda. \tag 6.6\endalign$$
The
last term can be ignored, since it does not make influence to
$\mu_U^{x_{\Lambda\backslash i}}$, and so neither
$\lambda_1^{x_{\Lambda\backslash i}}$. The coefficient of the second to
the last term varies over whole $\Bbb R$ if $J\ne 0$.

We consider two models only: $u(x)=\alpha x^2$ and $u(x)=
x^4-\beta x^2$ for some constants $\alpha>0$ and $\beta \in \Bbb
R$, respectively.

\proclaim{Theorem 6.2} Let $u(x)=\alpha x^2$ for some constant
$\alpha >0$ and let $U(x)=\sum_i u(x_i)+H(x)$ with Hamiltonian $H(x)=J\sum_{\langle ij\rangle }  (x_i-
x_j)^2$. Then we have
$$\inf_{\Lambda\Subset
{\Bbb Z}^d} \inf_{\omega\in {\Bbb R}^{\Bbb Z^d}} \lambda_1
\big(U_{\Lambda }^\omega \big)\ge \inf_{\Lambda\Subset
{\Bbb Z}^d}\inf_{\omega\in {\Bbb R}^{\Bbb Z^d}} \sz \big(U_{\Lambda
}^\omega \big) \ge 2\alpha. \tag 6.7$$
\endproclaim

\demo{Proof} It suffices to prove the second estimate. By Example 5.2 and Theorem 1.3, the
proof is very much the same as proving
$$\inf_{\Lambda\Subset
{\Bbb Z}^d} \inf_{\omega\in {\Bbb R}^{\Bbb Z^d}} \lambda_1
\big(U_{\Lambda }^\omega \big) \ge 2\alpha.$$
Hence we prove here the last assertion only. First,
we have
$$|\partial_{ij}U (x)| =
\cases
  2 J, &\qquad i, j\in \Lambda , \; |i-j|=1\\
  0, &\qquad i, j\in \Lambda , \; |i-j|> 1.
\endcases \tag 6.8$$ The right-hand side is independent of $x$, which is the
main reason why we were looking for the uniform estimates (with
respect to the linear term) in Examples 4.2 and 4.3. By (6.5), we
have $ \lambda_1^{x_{\Lambda\backslash i},\,\omega} \ge 2\alpha + 4 d J,$
which is indeed sharp in view of Example 4.2. Combining these
facts together and using (1.4) with $w_i\equiv 1$, it follows that
$$\align
\lambda_1 \big(U_{\Lambda }^\omega \big)& \ge \inf_{x\in \Bbb
R^{|\Lambda |}}\,\min_{i\in \Lambda } \bigg[2\alpha +4dJ -
\sum_{j\in \Lambda:\, |i-j|=1} 2 J\bigg]\\
& = 2\alpha
+4dJ - 2 J \max_{i\in \Lambda }\big|\{ \langle i,\, j\rangle:\,
j\in \Lambda\}\big|\\
& \ge 2\alpha \endalign$$uniformly in $\omega\in \Bbb
R^{\Bbb Z^d}$ and $\Lambda\Subset {\Bbb Z}^d$. The sign of the last
equality holds once $\llz$ contains a point together with all of its
neighbors.\qed\enddemo

In the last step of the proof, we did not use Theorem 1.2
since the matrix $\big(|\partial_{ij}U (x)|: i,j\in \llz\big)$
is very simple. Nevertheless, it provides us a good chance to justify
the power of Theorem 1.2. To do so, take $\eta_i^{x_{\Lambda\backslash i}}=2\az + 4 d J=\lambda_1^{x_{\Lambda\backslash i}}$. Then
$$\align
&s_i(x)=\eta_i^{x_{\Lambda\backslash i}}- \sum_{j\in\llz: j\ne i}|\partial_{ij}
U(x)|= 2\az + 4 d J- 2 J \big|\{\langle i, j\rangle: j\in \llz \}\big|,\qqd i\in \llz,\\
&{\underline s} (x)=\min_{i\in\llz} s_i(x)= 2\az.\endalign$$
Since $h^{(\gamma) }\ge 0$, Theorem 1.2 already gives us
$\lambda_1 \big(U_{\Lambda }^\omega \big)\ge \inf_{x}{\underline s} (x)=2\az $ as expected,
without using $h^{(\gamma) }$. To see the role played by $h^{(\gamma) }$, note that
$$\align
&q_i(x)=\eta_i^{x_{\Lambda\backslash i}}-{\underline s}(x)=4 d J,\qqd i\in \llz,\\
&d_i(x)= s_i(x)- {\underline s}(x)=4dJ  - 2 J \big|\{\langle i, j\rangle: j\in \llz \}\big|,\qqd i\in \llz.
\endalign$$
Note that $d_i(x)$ here depends on $i$. Thus
$$\align
h^{(\gamma) }(x)&=\min_{A:\, \emptyset \ne A\subset \llz
}\frac{1}{|A|} \bigg[\sum_{i\in A} \frac{d_i(x)}{q_i(x)^\gamma }+
\sum_{i\in A,\, j\in\llz\setminus A} \frac{|\partial_{ij} U
(x)|}{[q_i(x)\vee q_j(x)]^\gamma}\bigg]\\
&=\frac{2J}{(4dJ)^\gz}\min_{A:\, \emptyset \ne A\subset \llz
}\frac{1}{|A|}\sum_{i\in A}\Big[2d-\big|\{\langle i, j\rangle: j\in
\llz \}\big|
  +\big|\{\langle i, j\rangle: j\in \llz\setminus A \}\big|\Big]\\
 &=\frac{2J}{(4dJ)^\gz}\min_{A:\, \emptyset \ne A\subset \llz
}\frac{1}{|A|}\sum_{i\in A}\Big[
  \big|\{\langle i, j\rangle\}\big|-\big|\{\langle i, j\rangle: j\in A \}\big|\Big]\\
&=\frac{2J}{(4dJ)^\gz}\min_{A:\, \emptyset \ne A\subset \llz }\frac{1}{|A|}\sum_{i\in A}
  \big|\{\langle i, j\rangle: j\notin A \}\big|.\\
  &=:\frac{2J}{(4dJ)^\gz}\min_{A:\, \emptyset \ne A\subset \llz }\frac{|\partial A|}{|A|}.
\endalign  $$
Clearly, the right-hand side depends reasonably on the geometry of $\llz$.
Roughly speaking, by the isoperimetric principle, the last minimum of the ratio is approximately
${|\partial B|}/{|B|}$, where $B$ is the largest ball contained in $\llz$. Anyhow,
for regular $\llz$ (cube for instance),
$$h^{(\gamma) }(x)\le
\frac{2J}{(4dJ)^\gz }\cdot \frac{|\partial\llz|}{|\llz|}\to 0 \qqd\text{as } \llz\uparrow {\Bbb Z^d}.$$
Hence for this model,
 $h^{(\gamma) }$ makes no
contribution to $\lambda_1 \big(U_{\Lambda }^\omega \big)$ for the estimate uniformly in $\llz$.

\proclaim{Theorem 6.3} Let $u(x)=x^4-\beta x^2$ for some constant
$\beta\in \Bbb R$ and let $U(x)=\sum_i u(x_i)+H(x)$ with Hamiltonian $H(x)=-2 J\sum_{\langle
ij\rangle }  x_i  x_j$. Then we have
$$\align
{\hskip-2.5truecm}\inf_{\Lambda\Subset {\Bbb Z}^d}\inf_{\omega\in
{\Bbb R}^{\Bbb Z^d}} \lambda_1 \big(U_{\Lambda }^\omega \big)
&\ge\inf_{\Lambda\Subset {\Bbb Z}^d}\inf_{\omega\in {\Bbb R}^{\Bbb
Z^d}} \sz \big(U_{\Lambda
}^\omega \big)\\
& \ge \frac{\sqrt{\beta^2+8}- \beta}{\sqrt{e}}
\exp\bigg[-\frac{1}{8}\beta\Big(
\beta+\sqrt{\beta^2+8}\Big)\bigg]-4dJ, {\hskip-1truecm}\tag
6.9\endalign$$ For simplicity, we write $r=2dJ$. The right-hand side
is positive if $(\beta, r)\in \Bbb R\times \Bbb R_+$ is located in
the region below the curve in Figure 1 (including the region of
$\bz\le 0$ vertically below the shade one.)
\endproclaim

\demo{Proof}
As shown in part (2) of
Example 4.6, for zero boundary condition $\omega=0$, we have
$$\lim_{\beta\to\infty}\sz^{x_{\Lambda\backslash i},\,\omega}\le \lim_{\beta\to\infty}\lambda_1^{x_{\Lambda\backslash i},\,\omega}=0.$$
 In
other words, due to the double-well potential, the spectral gap and then the logarithmic constant
will be absorbed as $\beta\to\infty$. Combining Example 5.3 with Theorem 1.3
and following the last step of the proof
Theorem 6.2, we obtain the required lower estimate. \qed
 \input epsf.sty
$$\epsfbox[140 600 390 770]{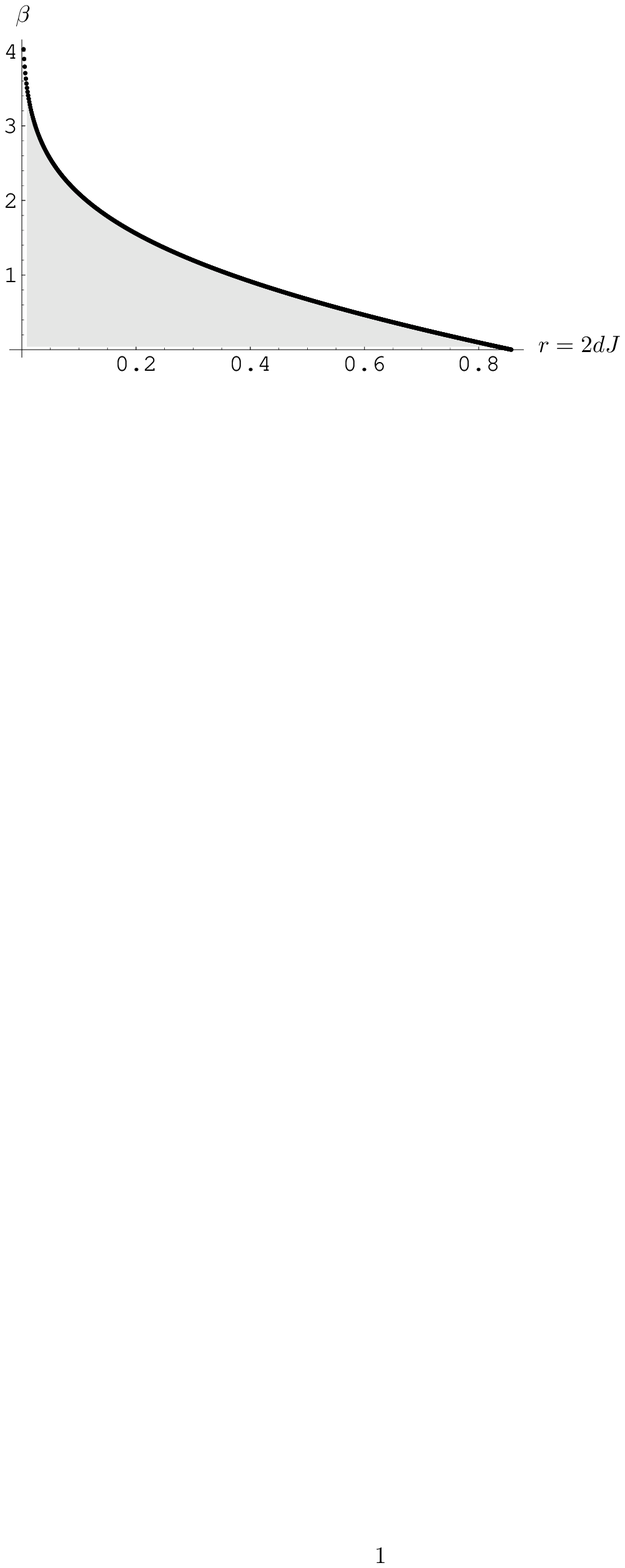}$$
\centerline{Figure 1.}
\enddemo

For the Hamiltonian $H(x)=J\sum_{\langle
ij\rangle }  (x_i - x_j)^2$ discussed several times before,
simply replacing $\bz$ with $\bz- 2d J$ in Theorem 6.3, we obtain the following
estimate:
$$\align
\inf_{\Lambda\Subset
{\Bbb Z}^d}\,&\inf_{\omega\in {\Bbb R}^{\Bbb Z^d}} \lambda_1
\big(U_{\Lambda }^\omega \big)\\
&\ge\inf_{\Lambda\Subset
{\Bbb Z}^d}\inf_{\omega\in {\Bbb R}^{\Bbb Z^d}} \sz \big(U_{\Lambda
}^\omega \big)\\
& \ge \frac{\sqrt{(\beta-r)^2+8}- \beta+r}{\sqrt{e}}
\exp\bigg[-\frac{1}{8}(\beta-r)\Big( \beta-r+\sqrt{(\beta-r)^2+8}\Big)\bigg]\\
&\qd -2r, \tag 6.10\endalign$$ where $r=2 d J$. The ergodic region
is shown in Figure 2.

$$ \epsfbox[125 605 390 780]{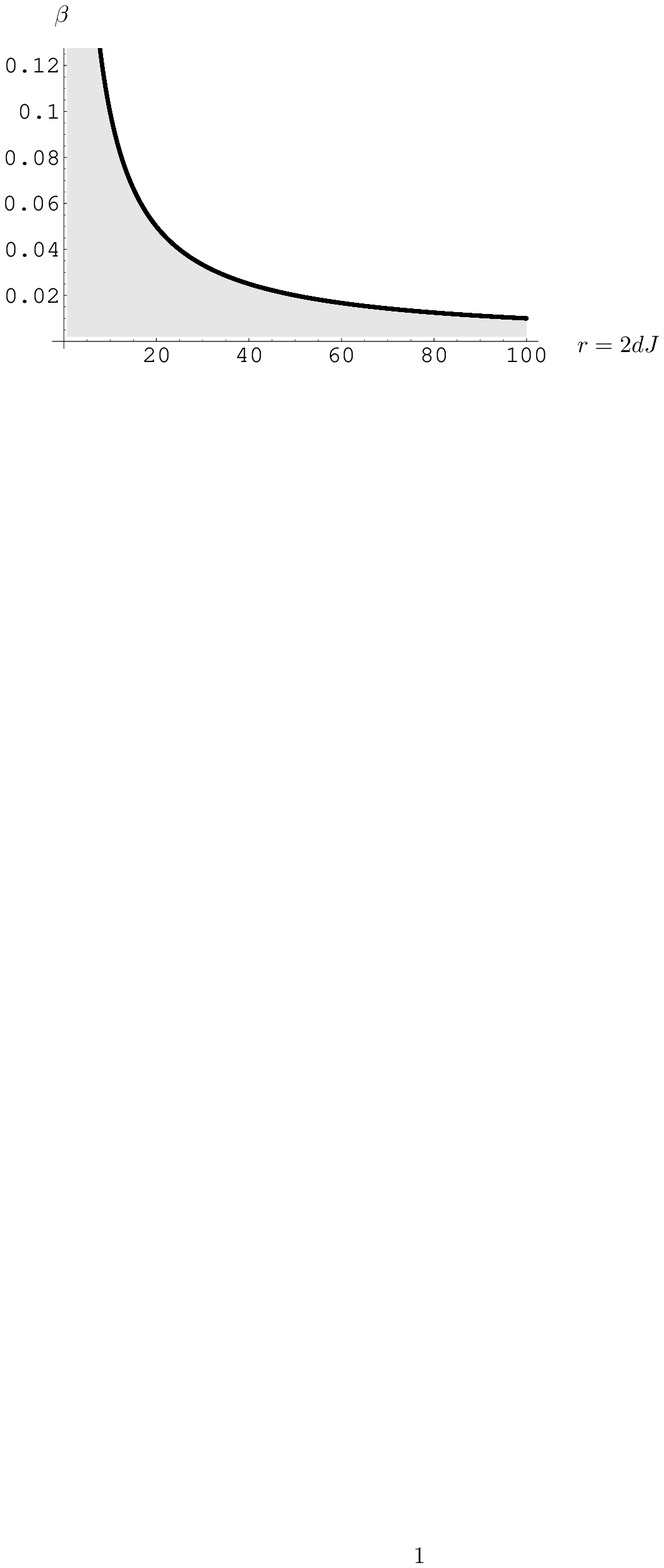}$$
\centerline{Figure 2.}

\proclaim{Remark 6.4} {\rm As mentioned below the proof of Proposition 3.1, by considering the interacting
terms more carefully, one may improve Theorem 1.1 for stronger interactions. For instance, since
the variance of a random variable having the distribution with density $\exp[-x^4+\beta x^2]/Z$ on the real line
is asymptotically $\beta/2$ for $\beta\ge 0$, and is bounded above by
$$\frac{\ggz(3/4)}{\ggz(1/4)}+{\bz}\bigg/\bigg[{2+\frac{4\ggz(1/4)}{9(1+\bz)\ggz(3/4)}}\bigg],$$
by using [9; Proposition 5.8], when $\beta\ge 0$, the lower bound of
$\inf\limits_{\Lambda\Subset {\Bbb Z}^d}\inf\limits_{\omega\in {\Bbb
R}^{\Bbb Z^d}} \lambda_1 \big(U_{\Lambda }^\omega \big)$ given in
Theorem 6.3 can be improved as follows: replacing the interaction
term $4dJ$ in $(6.9)$ with
$$\align
&4dJ\bigg[\frac{\ggz(3/4)}{\ggz(1/4)}\!+\!{\bz}\bigg/\bigg[{2\!+\!\frac{4\ggz(1/4)}{9(1\!+\!\bz)\ggz(3/4)}}\bigg]\bigg]
\frac{\sqrt{\beta^2+8}\!-\! \beta}{\sqrt{e}}
\exp\bigg[\!-\!\frac{1}{8}\beta\Big( \beta\!+\!\sqrt{\beta^2\!+\!8}\Big)\bigg].\\
&\tag 6.11\endalign$$}
\endproclaim

Finally, we mention that there is another technique which works even
in the irreversible situation (cf. [17; Theorem 14.10]) to handle
with the exponentially ergodic region, because the second model
(Theorem 6.3) is attractive (stochastic monotone) and has the
moments of all orders, plus a use of the translation invariant.
However, as known that the logarithmic Sobolev inequality already
implies an exponential ergodicity in the entropy and moreover, the
usual exponential ergodicity is equivalent to the Poincar\'e
inequality with nearly the same convergence exponent in the present
context (cf. [13; Theorem 8.13]), there is almost no room to improve
the ergodic region.
\bigskip

\flushpar{\bf Acknowledgements}. This paper takes an unusual long period in preparation.
The most part of the paper was finished in 2002, but the exact
coefficient $1/4$ of decay rate $\exp[-\bz^2/4]$ for the second model (Theorem 6.3) was left to be open and so the earlier draft was
communicated within a small group only. Recently, the author came back to compute the logarithmic Sobolev
constant which leads to the precise coefficient and hence completes the paper.
The author would like to acknowledge the organizers for several conferences in which partial results of the paper
were presented: I. Shigekawa (the 11th International Research Institute of Mathematical Society of Japan, 2002),
L.M. Wu (Chinese and French Workshop on Probability and Applications, 2004), Z.M. Ma and M. R\"ockner (Second
Sino-German Meeting on Stochastic Analysis, 2007). The author is also greatly indebted to
the referees for their helpful comments.


\bigskip

\flushpar Note: Figures 1 and 2 were missed in
the publication and the correction appeared in
the same journal, Vol. 25, No. 12, pp. 2199-2199.

\enddocument